\documentclass{article}

\usepackage{times}
\usepackage[left=3cm,top=2cm,right=3cm,bottom=2cm]{geometry}
\usepackage{amsmath}
\usepackage{enumerate}
\usepackage{amssymb}

\usepackage{amsfonts}
\usepackage[small,nohug,heads=littlevee]{diagrams}
\usepackage[all]{xy}

\usepackage{amsthm}

\theoremstyle{definition}
\newtheorem{definition}{Definition}[section]
\theoremstyle{plain}
\newtheorem{theorem}[definition]{Theorem}
\newtheorem{lemma}[definition]{Lemma}
\newtheorem{proposition}[definition]{Proposition} 
\newtheorem{corollary}[definition]{Corollary}

\theoremstyle{remark}
\newtheorem{remark}[definition]{Remark} 
\newtheorem{example}[definition]{Example}

\newcommand{\rC}{\mathbb C} 
\newcommand{\bP}{\mathbb{P}} 
\newcommand{\sO}{\mathcal O} 
\newcommand{\vp}{\varphi} 
\newcommand{\tC}{\tilde{C}}

\newcommand{\W}{\Omega}
\newcommand{\tX}{\tilde{X}}
\newcommand{\tY}{\tilde{Y}}
\newcommand{\tM}{\tilde{M}}

\newcommand{\ev}{\text{ev}}
\newcommand{\Monxb}{\overline{\mathcal{M}_{0,n}}(X,\beta)}
\newcommand{\Mgnxb}{\overline{\mathcal{M}_{g,n}}(X,\beta)}
\newcommand{\Mvir}{\left[ \Mgnxb \right]^{\mathrm{vir}}}
\newcommand{\vir}{\mathrm{vir}}
\newcommand{\codim}{\mathrm{codim}_{\rC}}
\newcommand{\rk}{\mathrm{rk}}
\newcommand{\cD}{\mathcal{D}}
\newcommand{\deri}{\mathcal{D}}
\newcommand{\Edot}{\mathcal{E}_{\bullet}}
\newcommand{\Fdot}{\mathcal{F}_{\bullet}}
\newcommand{\E}{\mathcal{E}}
\newcommand{\F}{\mathcal{F}}
\newcommand{\C}{\mathcal{C}}
\newcommand{\tp}{\tilde{p}}
\newcommand{\tf}{\tilde{f}}
\newcommand{\sL}{\mathcal{L}}
\newcommand{\Ldot}{\mathcal{L}_{\bullet}}

\newcommand{\truntwo}{\tau_{\leq 2}}
\newcommand{\Lup}{\mathcal{L}^{\bullet}}
\newcommand{\Eup}{\mathcal{E}^{\bullet}}
\newcommand{\Fup}{\mathcal{F}^{\bullet}}

\newcommand{\Dptwo}{\cD(p)_2}
\newcommand{\Dptwoset}{ \{ \Dptwo \}_{p\in \bM}}
\newcommand{\vect}{\text{Vect}}

\newcommand{\vectLvF}{\vect(\sL_2 \oplus \vp^* \F_2)}
\newcommand{\vectvF}{\vect (\vp^* \F_2)}
\newcommand{\vectF}{\vect(\F_2)}
\newcommand{\zcone}{0_{\sL_2}^![\C_{\bM}]}

\newcommand{\supp}{\text{supp}}

\newcommand{\vdim}{\mathrm{vdim}_{\rC}}
\newcommand{\bM}{\mathbf{M}}
\newcommand{\bN}{\mathbf{N}}
\newcommand{\bMN}{[\bM, \bN]^\vir}
\newcommand{\sC}{\mathfrak{C}}
\newcommand{\moduli}{\overline{\mathcal{M}}}
\newcommand{\tbeta}{\tilde{\beta}}
\newcommand{\tGamma}{\tilde{\Gamma}}
\newcommand{\univ}{\text{univ}}
\newcommand{\finite}{\text{finite}}
\newcommand{\red}{\text{red}}
\newcommand{\calA}{\mathcal{A}}
\newcommand{\calU}{\mathcal{U}}
\newcommand{\calW}{\mathcal{W}}
\newcommand{\calF}{\mathcal{F}}
\newcommand{\comp}{\text{comp}}
\newcommand{\coarse}{\overline{M}}
\newcommand{\longvec}[1]{\stackrel\longrightarrow{\smash{#1}\vphantom{i}\,}}

\title{Gromov-Witten invariants of blow-ups along submanifolds with convex normal bundles}
\author{Hsin-Hong Lai}
\date{}

\begin{document}

\maketitle

\begin{abstract}
When the normal bundle $N_{Z/X}$ is convex with a minor assumption, we prove that genus$-0$ GW-invariants of the blow-up $Bl_Z X$ of $X$ along a submanifold $Z$, with cohomology insertions from $X$, are identical to GW-invariants of $X$. Under the same hypothesis, a vanishing theorem is also proved. An example to which these two theorems apply is when $N_{Z/X}$ is generated by its global sections.  These two main theorems do not hold for arbitrary blow-ups, and counter-examples are included.
\end{abstract}

\tableofcontents

\section{Introduction}

In~\cite{R1}, Y. Ruan proposes naturality problems of quantum cohomology rings under birational surgery. In~\cite{Ruanuniruled}~\cite{McDuff}, GW-invariants are used to classify symplectic manifolds in a symplectic birational geometric program. Recently, there has also been substantial progress in crepant resolution conjecture. On the other hand, blow-up formula for GW-invariants is known only for very few cases. Let $\pi: \tX \to X$ be the blow up of $X$ along the submanifold $Z$. A natural question is if the induced genus$-0$ GW-invariants of $\tX$ coincide with the GW-invariants of $X$. That is, if $\alpha_i \in H^*(X)$ and $\beta\in H_2(X)$, do we have
\begin{equation}{\label{Q1}}
\langle \pi^*\alpha_1, \cdots, \pi^*\alpha_n \rangle _{0,n, \pi^!\beta}^{\tX} = \langle \alpha_1, \cdots, \alpha_n \rangle _{0,n, \beta}^X ?
\end{equation}
When formulated in this generality, the answer is negative (see Remark 9 in~\cite{vertex} or Example~\ref{counter1}).
In~\cite{Gathmann}, \cite{Hu} and \cite{Hu1}, the answer to Question~\eqref{Q1} has been shown to be true in some cases, where $\dim Z \leq 2$ with various assumptions, including the requirement that cohomology insertions are supported away from $Z$ when $\dim Z =2$.
In this paper, we will show that if the normal bundle $N_{Z/X}$ is convex with a minor assumption, then the answer to Question~\eqref{Q1} is also affirmative. This provides examples where $\dim Z$ can be any number without assuming cohomology insertions are supported away from $Z$. First recall the definition of a convex bundle:
\begin{definition}
A vector bundle $W$ over a manifold $Z$ is called convex if and only if $H^1(\mathbb{P}^1, f^*W)=0$ for any holomorphic map $f: \mathbb{P}^1 \to Z$.
\end{definition}
In this paper, we consider two classes of submanifolds $Z\subset X$.
\begin{definition}{\label{0typeI}}
A connected submanifold $Z \subset X$ is of type I, if the following two conditions are satisfied:
\begin{enumerate}
\item $N_{Z/X}$ is a convex bundle over Z,
\item There is a subbundle $\calF$ in $N_{Z/X}$ with rank $\rk (\calF)\geq 2$, and $\calF$ is generated by global sections.
\end{enumerate}
\end{definition}
An example of type I is when $N_{Z/X}$ is generated by global sections.
\begin{definition}{\label{0typeII}}
A connected submanifold $Z \subset X$ is of type II, if every holomorphic map $f : \bP^1 \to Z$ must be a constant map.
\end{definition}
For example, $Z$ is of type II if $Z$ is a product of higher genus curves or abelian varieties.

Our first main result is the following:
\begin{theorem}{\label{thm1}}
Suppose each connected component of the submanifold $Z = \coprod_{i} Z_i \subset X$ is of type I or type II.
Let $V$ be a vector bundle over $X$, and $\mathbf{c}$ be an invertible multiplicative characteristic class. Then we have an equality of genus-$0$ twisted Gromov-Witten invariants
$$\langle \alpha_1, \cdots, \alpha_n \rangle _{0,n, \beta}^{X, \mathbf{c}, V}=\langle \pi^*\alpha_1, \cdots, \pi^*\alpha_n \rangle _{0,n, \pi^!\beta}^{\tX, \mathbf{c}, \pi^*V}, \text{ where } \alpha_i \in H^*(X) \text{ for all } i.$$
\end{theorem}

Given an arbitrary projective manifold $X$, Example~\ref{source} provides several ways to find a submanifold $Z\subset X$, so that $N_{Z/X}$ is generated by global sections. This is the major source of examples to which Theorem~\ref{thm1} applies. Type I and type II cases cover most cases when $N_{Z/X}$ is convex. We speculate that Theorem~\ref{thm1} holds as long as $N_{Z/X}$ is convex without any additional assumptions. 

Convexity of the normal bundle is a critical assumption in Theorem~\ref{thm1}. This is illustrated
by Example~\ref{counter1}, which has the following properties:\\
(1) The submanifold $Z\subset X$ has enough freedom to move inside $X$, so that $Z$ can avoid any finite collection of holomorphic curves.\\
(2) The moduli spaces of $\tX$ and $X$ are both smooth and birational to each other.\\
(3) The difference of (push-down) virtual classes has non-zero contribution to GW-invariants. Therefore the conclusion of Theorem~\ref{thm1} does not hold in this case.\\
In this example, the non-convex part of the normal bundle $N_{Z/X}$ "twists" the obstruction bundle on the moduli space of $\tX$, and gives rise to the correction term of (push-down) virtual classes/GW-invariants.

Theorem~\ref{thm1} is a direct consequence of the following equality of virtual classes. $\tilde{W_0}$ and $W_0$ are degenerations (from deformation to the normal cones) of $\tX$ and $X$ respectively.

\begin{theorem}{\label{thm2}}
Suppose each connected component of the submanifold $Z = \coprod_{i} Z_i \subset X$ is of type I or type II. Then we have $\phi_*[\moduli (\tilde{\mathcal{W}_0}, 0, n , \pi^! \beta)]^\vir 
= [\moduli (\mathcal{W}_0, 0, n , \beta)]^\vir$.
\end{theorem}

In some special cases, Theorem~\ref{thm2} can be improved as follows:

\begin{theorem}{\label{trans0}}
Suppose $Z$ is the transversal intersection of two arbitrary manifolds $X$ and $Y$ in a compact homogeneous space $\mathcal{P}$.
Then we have $\vp_*[\overline{\mathcal{M}_{0,n}}(Bl_Z X, \pi^!\beta)]^{\vir} = [\Monxb]^{\vir}$ in the Chow group.
\end{theorem}
As a corollary, if $X$ is an arbitrary projective manifold and $Z$ is a collection of points, then the equality of virtual classes holds. The case where $X$ is a convex manifold and $Z$ is a collection of points, has been proved in~\cite{Gathmann}. We remark that when $g>0$ and $Z$ is a point, in general we have $\vp_*[\overline{\mathcal{M}_{g,n}}(Bl_Z X, \pi^!\beta)]^{\vir} \neq [\overline{\mathcal{M}_{g,n}}(X, \beta)]^{\vir}$.

The second part of this paper is a vanishing theorem. First we introduce some notation.
\begin{enumerate}
\item[$\bullet$] $[n]:= \{ 1,2,\cdots,n\}.$
\item[$\bullet$] Given $A \subset [n]$, use $\longvec{\tau_\bullet \alpha_A}$ to denote descendant insertions $\{ \tau_{i_a}\cdot \alpha_a\}_{a\in A}$, where $\alpha_a\in H^*(X)$ and $i_a \geq 0$. If $i_a=0$ for all $a\in A$, then $\longvec{\tau_\bullet \alpha_A}$ is simply denoted by $\longvec{\alpha_A}$.
\item[$\bullet$] $\longvec{\mathbf{1}_{[n]}}:=(1,1,\cdots,1)$, where $1\in H^*(X)$.
\item[$\bullet$] The product 
$\longvec{\tau_\bullet \alpha_A} \cdot \longvec{\tau_\bullet \gamma_B} := \{ \tau_{i_a+j_b}\cdot \alpha_a \cap \gamma_b   \}_{a=b\in A\cap B} \cup \{ \tau_{i_a}\cdot \alpha_a\}_{a\in A-B} \cup \{ \tau_{j_b}\cdot \gamma_b\}_{b\in B-A}.$
\item[$\bullet$] The GW-invariant $\langle \longvec{\tau_\bullet \alpha_A} \cdot \longvec{\mathbf{1}_{[n]}} \rangle_{0,n,\beta}^X$ is simply denoted by $\langle  \longvec{\tau_\bullet \alpha_A} \rangle_{0,n, \beta} ^X$.
\end{enumerate}
\begin{theorem}{\label{thm3}}
$I, J, K$ are disjoint sets with $J \subset [n]$.
Suppose $Z = (\coprod_{i\in I} Z_i) \cup (\coprod_{j\in J} Z_j) \cup (\coprod_{k\in K} Z_k)$ is a disjoint union of submanifolds in $X$, with the following assumptions:
\begin{enumerate}
\item[$\bullet$] For each $i \in I\cup J$, $Z_i \subset X$ is either of type I or of type II.
\item[$\bullet$] For each $k \in K$, $N_{Z_k/X}$ is convex.
\item[$\bullet$] The curve class $\tbeta = \pi^!\beta + \sum_{i\in I} d_i e_i + \sum_{j\in J} d_j e_j +\sum_{k\in K} d_k e_k$ with $d_i\neq 0$ for all $i\in I$, and $0 \neq \beta\in H_2(X)$. Here $e_\bullet$ are the exceptional line classes.
\item[$\bullet$] $\longvec{\omega_J}$ is a collection of cohomology classes in $H^*(\tX)$. And $PD_{\tX}(\omega_j)$ lies in the image of $H_*(E_j) \to H_*(\tX)$, where $E_j$ is the exceptional divisor. 
\end{enumerate}
For $i\in I\cup J$, define 
$$
\delta_i= \left\{
\begin{array}{ll}
\rk(\calF) -1 & \text{, if } Z_i \subset \text{ X is of type I, and } \calF\subset N_{Z_i/X} \text{ is generated by global sections.}\\
\rk(N_{Z_i/X})-1 & \text{, if } Z_i \subset \text{ X is of type II.}
\end{array}
\right.
$$
Then 
$$\langle \longvec{\pi^*\alpha_A} \cdot \longvec{\tau_\bullet \gamma_{[n]}} \cdot \longvec{\omega_J} \rangle _{0, n, \tbeta}^{\tX}=0 \text{ when } \deg \longvec{\alpha_A} > 2\vdim \moduli_{0, A} (X, \beta)- 2\sum_{i\in I} \delta_i - 2\sum_{j\in J} \delta_j.$$
Here $\longvec{\alpha_A}$ is a collection of cohomology classes from $X$ with $A \subset [n]$, and $\longvec{\tau_\bullet \gamma_{[n]}}$ are arbitrary descendant insertions of $\tX$.
\end{theorem}

Roughly speaking, when taking $J=\emptyset$, Theorem~\ref{thm3} can be numerically interpretated as: $$\text{ The image of }\vp: \moduli_{0,n}(\tX, \tbeta) \to \moduli_{0,A}(X, \beta) \text{ has "virtual codimension" }\geq \sum_{i\in I} \delta_i.$$
Therefore, if there are too many cohomology insertions from $X$, then the GW-invariant of $\tX$ vanishes. In~\cite{GaThesis}, Gathmann proved a vanishing theorem for genus-0 non-descendant GW-invariants when blowing up at points.
Theorem~\ref{thm3} is a generalization of Gathmann's results in two aspects:
$$
\begin{array}{ll}
(1)\text{ There is no restriction on } \dim Z. & (2) \text{Theorem~\ref{thm3} also holds for descendant GW-invariants. }
\end{array}
$$
We remark that Theorem~\ref{thm3} only holds for blow-ups with convex normal bundles, but does not hold for arbitrary blow-ups (see Example~\ref{counter2}).

In Example~\ref{surface}, we use Theorem~\ref{thm3} to show that, given any algebraic surface $S$ which is not (birationally equivalent to) a ruled or rational surface, then most genus$-0$ descendant GW-invariants of $S$ are zero. When $p_g(S)>0$, this conclusion has been deduced from the Image Localization Theorem of holomorphic two forms in~\cite{structure}.

The tools used in this paper are : degeneration formula (~\cite{Ruandeg}~\cite{IP}~\cite{deg}~\cite{topGW}), compatibility of perfect obstruction theories (see Definition~\ref{defcompatible} and~\cite{inc}\cite{functor}\cite{LT2}) and deformation invariance of virtual classes. Since there is no assumption on the manifold $X$, the moduli of stable maps of $X$ can be highly singular. Instead of analyzing singularities of the moduli space (which is nearly impossible), in Section 3 we show that if $N_{Z/X}$ is convex, then $\moduli_{0,n}(\tX, \tbeta) \to \moduli_{0,n}(X, \pi_*\tbeta)$ has compatible perfect obstruction theories. General blow-ups don't have this property. We use Proposition~\ref{special} as a criterion to the equality of (push-forward) virtual classes.

To prove Theorem~\ref{trans0}, we deform the submanifold $Z$ so that the technical assumption in Proposition~\ref{special} is satisfied. Regarding the type I case in Theorem~\ref{thm2}, degeneration formula (in cycle forms) is used to split the problem into various relative virtual classes associated to a ruled variety $\bP_Z(N_{Z/X}\oplus \sO_Z)$, and then the submanifold $Z$ is moved so that the technical assumption in Proposition~\ref{special} is satisfied. For type II case in Theorem~\ref{thm2}, we move holomorphic curves instead of $Z$ and argue directly. Although one can always move holomorphic curves as long as $N_{Z/X}$ is convex, there is a technical difficulty in applying Proposition~\ref{special} due to singularities of the moduli space. See Remark~\ref{difficulty1} for discussion.

Our starting point for the vanishing theorem is Lemma~\ref{vanish}, which also requires compatible perfect obstruction theories, and therefore doesn't hold for arbitrary blow-ups. The bound of the degree of cohomology insertions in Theorem~\ref{thm3}, is deduced from codimension analysis of the image on virtual normal cones.\\

When $N_{Z/X}$ is a direct sum of convex and concave bundles, in general we have
$$\langle \pi^*\alpha_1, \cdots, \pi^*\alpha_n \rangle _{0,n, \pi^!\beta}^{\tX} \neq \langle \alpha_1, \cdots, \alpha_n \rangle _{0,n, \beta}^X.$$
The correction term will be discussed in the future.

\section*{Acknowledgements}
I would like to thank Harry Tamvakis for teaching me Gromov-Witten theory, Bong Lian for kindly suggesting me to consider virtual classes of blow-up at points, and Daniel Ruberman for his support during the course of this work. I would also like to thank Dan Abramovich and Jun Li for helpful conversations.

\section{Preliminaries and notation}
Given a projective manifold $X$ and a curve class $\beta \in H_2(X)$, the stable maps moduli $\moduli_{g,n}(X, \beta)$ collects all holomorphic map from a genus-$g$ nodal curve with $n$ marked points $f: C \to X$. These holomorphic maps are required to satisfy the stability condition, which means the automorphism of each map is finite. Let $\mathcal{C}:= \moduli _{g,n+1}(X, \beta)$ be the universal curve of $M:= \moduli _{g,n}(X, \beta)$. Recall that the perfect tangent obstruction complex of $\moduli _{g,n}(X, \beta)$ is given by 
$$\Fdot = [\F_1 \to \F_2] = [\E xt_{\C/M}^{\bullet} ([f^*\W_X \to \W_{\C/M}(D)],\sO_{\C})],$$ where $f: \C \to X$ is the universal map and $D$ are the marked sections of $\moduli _{g,n}(X, \beta)$.
One also has:
\begin{enumerate}
\item a evaluation map $\text{ev}: \moduli_{g,n}(X, \beta) \to X^n$, which evaluates at the marked points,
\item a line bundle $\mathbb{L}_i$ with the fiber over $(C,a_1,\cdots,a_n,f)$ isomorphic to the cotangent space of $C$ at $a_i$.
\end{enumerate}

Let $\psi_i$ be the first Chern class $c_1(\mathbb{L}_i)$.
Given $\gamma_i \in H^*(X)$, for $i=1,\cdots,n$, the genus-$g$ descendant Gromov-Witten invariants are defined as:
$$ \langle \tau_{a_1}\gamma_1, \cdots, \tau_{a_n}\gamma_n \rangle _{g,n, \beta}^X= \int_{\Mvir}
\psi_1^{a_1}\cap \cdots \cap \psi_n^{a_n}\cap \ev^*(\otimes_{i=1}^n \gamma_i).$$

Suppose $V$ is a vector bundle over $X$. Consider the universal family:
$$
\xymatrix{
\overline{\mathcal{M}_{g,n+1}}(X,\beta) \ar[r]^>>>{e_{n+1}} \ar[d]^{\pi_{n+1}}& X\\
\overline{\mathcal{M}_{g,n}}(X, \beta)
}
$$
$(R\pi_{n+1})_* \circ e_{n+1}^*(V)$ can be represented by a two-term complex of vector bundles $[V_0 \to V_1]$. If $\mathbf{c}$ is an invertible multiplicative characteristic class, the twisted genus-$g$ descendant Gromov-Witten invariants defined in~\cite{Twisted} are given by:
$$ \langle \tau_{a_1}\gamma_1, \cdots, \tau_{a_n}\gamma_n \rangle _{g,n, \beta}^{X, \mathbf{c}, V}= \int_{\Mvir}
\psi_1^{a_1}\cap \cdots \cap \psi_n^{a_n}\cap \ev^*(\otimes_{i=1}^n \gamma_i) \cap \mathbf{c}(V_0 \ominus V_1).$$

\section{Blow-ups with convex normal bundles}
\subsection{Compatibility of perfect obstruction theories}
Given any morphism $\pi : Y \to X$ of two projective manifolds and $\tbeta \in H_2(Y)$, there always exists an induced map $\vp : \moduli_{g,n}(Y, \tbeta) \to \moduli_{g,n}(X, \pi_*\tbeta)$, as long as $\moduli_{g,n}(X, \pi_*\tbeta)$ makes sense (this is equivalent to saying $n \geq 3$ if $\pi_*\tbeta =0$). Suppose $\Edot = [\E_1 \to \E_2]$ and $\Fdot = [\F_1 \to \F_2]$ are the perfect tangent-obstruction complexes on $\moduli_{g,n}(Y, \tbeta)$ and $\moduli_{g,n}(X, \pi_*\tbeta)$ respectively, there always exists a natural map $\Edot \to \vp^*\Fdot$ in $\mathcal{D}(\mathcal{O}_{\moduli_{g,n}(Y, \tbeta)})$, the derived category of the coherent sheaves on $\moduli_{g,n}(Y, \tbeta)$.
The obstruction sheaves of $\Edot$ on $\moduli_{g,n}(Y, \tbeta)$ and $\Fdot $ on $\moduli_{g,n}(X, \pi_*\tbeta)$ are defined as $\text{Ob}_{\moduli_{g,n}(Y, \tbeta)}:=h^2(\Edot)$ and $\text{Ob}_{\moduli_{g,n}(X, \pi_*\tbeta)}:=h^2(\Fdot)$. There is a natural map $\text{Ob}_{\moduli_{g,n}(Y, \tbeta)} \to \vp^*(\text{Ob}_{\moduli_{g,n}(X, \pi_*\tbeta)})$.

Suppose we have a stable map $\tilde{p}=(\tilde{C}, \underline{\tilde{a}}, \tilde{f}) \in \moduli_{g,n}(Y, \tbeta)$, where $\tilde{C}$ is a nodal curve and $\underline{\tilde{a}} \subset \tilde{C}$ are the marked points. The composition $\pi\circ\tilde{f}: (\tilde{C}, \underline{\tilde{a}}) \to X$ might not be stable. One contracts the unstable components to obtain the domain curve $C$. Then $\vp(\tilde{p}) \in \moduli_{g,n}(X, \pi_*\tbeta)$ is given by $(C, \underline{a}, f)$, where $\underline{a} \subset C$ are the marked points after contraction. We have the commutative diagram:
$$
\xymatrix{
\tC \ar[r]^{\tf} \ar[d]_{\psi=\text{stablization}}& Y \ar[d]^{\pi}\\
C \ar[r]_f &X
}
$$

\begin{lemma}{\label{can iso}}
There are canonical isomorphisms:

\begin{enumerate}
\item $H^0(C, f^*TX) \cong H^0(\tilde{C}, \tilde{f}^*\pi^*TX)$
\item $H^1(C, f^*TX) \cong H^1(\tilde{C}, \tilde{f}^*\pi^*TX)$.
\end{enumerate}

\end{lemma}

\begin{proof}
\hspace{.5in}

$$\tC \rTo^{\psi} C \rTo^{pt} \text{point}$$

Set $\mathcal{F} = f^*TX$, then there is a natural morphism $\mathcal{F} \rightarrow R\psi_*\circ L\psi^*(\mathcal{F})$. For any connected component $C_i$ of $C_\text{unstab}$, it must be a genus-$0$ nodal curve. Therefore

$$ H^0(C_i, \mathcal{O}_{C_i}) = \mathbb{C} \text{   and   } H^j(C_i, \mathcal{O}_{C_i}) =0 \text{  ,  for  } j\neq 0.$$
This implies $\mathcal{F} \rightarrow R\psi_*\circ L\psi^*(\mathcal{F})$ is an isomorphism.
Since $\tilde{C}$ and $C$ are proper, we have $\psi_! = \psi_*$ and $pt_! = pt_*$, therefore 

$$Rpt_*(\mathcal{F}) \rightarrow Rpt_*\circ R\psi_*\circ L\psi^*(\mathcal{F})= Rpt_!\circ R\psi_!\circ L\psi^*(\mathcal{F}) \cong R(pt\circ \psi)_!\circ L\psi^*(\mathcal{F}) $$
is an isomorphism. The lemma follows from cohomology of this isomorphism.
\end{proof}

Let $\pi : \tX :=Bl_Z X \to X$ be the blow-up of $X$ along $Z \subset X$. When $g=0$ and the normal bundle $N_{Z/X}$ is convex, we have surjectivity between obstruction sheaves.

\begin{proposition}{\label{surjective}}
If $N_{Z/X}$ is convex, then the natural map $\text{Ob}_{\moduli_{0,n}(\tX, \tbeta)} \to \vp^*(\text{Ob}_{\moduli_{0,n}(X, \pi_*\tbeta)})$ is surjective, where $\tbeta \in H_2(\tX)$.
\end{proposition}

\begin{proof}
For convenience, denote $\moduli_{0,n}(\tX, \tbeta)$ by $\tM$ and $\moduli_{0,n}(X, \pi_*\tbeta)$ by $M$.
Given a point $\tilde{p}=(\tilde{C}, \underline{\tilde{a}}, \tilde{f}) \in \tM$, the obstruction space is
$(\text{Ob}_{\tM})_{\tilde{p}} = h^2(\Edot \otimes_{\sO_{\tM}} k(\tp))$. We also have
$(Ob_M)_p = h^2(\Fdot \otimes_{\sO_M} k(p))$, where $p:=\vp(\tilde{p})=(C, \underline{a}, f) \in M$.
Consider the following commutative diagram of the right exact sequence (see Lemma~\ref{can iso}) :
$$
\xymatrix{
H^1(\tC, \tf^*T\tX) \ar[r] \ar[d]& (\text{Ob}_{\tM})_{\tp} \ar[dd]\ar[r] & 0\\
H^1(\tC, \tf^*\pi^*TX)\\
H^1(C, f^*TX)  \ar[u]_{\cong} \ar[r] & (\text{Ob}_{M})_{p} \ar[r] &0
}
$$
It suffices to prove 
$$H^1(\tilde{C}, \tilde{f}^*T\tX) \rightarrow H^1(\tilde{C}, \tilde{f}^*\pi^*TX) \text{ is surjective.}$$
First we pull back the blow-up exact sequence (see Lemma 15.4 in~\cite{Intersection})
$$0 \to T\tX \to \pi^* TX \to i_*Q \to 0$$
to $\tilde{C}$, where Q is the universal quotient bundle on the exceptional divisor $E=\bP (N_{Z/X})$:
$$\tf^*T\tX \to \tf^*\pi^*TX \to \tf^*Q \to 0.$$
And let $K_1$ and $K_2$ be the corresponding kernels
$$0 \to K_1 \to \tf^*\pi^*TX \to \tf^*Q \to 0$$
$$0 \to K_2 \to \tf^*T\tX \to K_1 \to 0$$
Since the domain curve has dimension$=1$, $H^2(\tC, K_2)=0$, which implies
$$H^1(\tC, \tf^*T\tX) \to H^1(\tC,K_1) \text{ is surjective.}$$
If we can show $H^1(\tC, \tf^*Q)=0$, then the composition $H^1(\tC, \tf^*T\tX) \to H^1(\tC, K_1) \to H^1(\tC, \tf^*\pi^*TX)$ is also surjective. Set $C':=\tf^{-1}(E)$, then $H^1(\tC, \tf^*Q) = H^1(C', \tf^*Q)$. Note $C'$ might be disconnected.
\\
There is another exact sequence on the exceptional divisor $E$ 
$$ 0 \to \sO_{N_{Z/X}}(-1) \to \pi^*(N_{Z/X}) \to Q \to 0, \text{ where } \pi: E \to Z.$$
Pull it back to $C'$ to deduce the right exact cohomology sequence
$$H^1(C', \tf^*\pi^*(N_{Z/X})) \to H^1(C',\tf^*Q) \to 0$$
Note that $C'$ is a collection of points and genus-$0$ nodal curves, and $N_{Z/X}$ is convex, therefore we have $H^1(C', \tf^*\pi^*(N_{Z/X}))= 0$. This implies $H^1(C', \tf^*Q)=0$ and completes the proof.

\end{proof}

In \cite{inc} and \cite{LT2}, the existence of global vector bundles is used to construct virtual fundamental classes. This technical assumption has been removed due to the work of A. Kresch ~\cite{Kr2}\cite{Kr1}. Nevertheless, for simplicity, in this paper we still assume the existence of global vector bundles, which is true in Gromov-Witten Theory (see~\cite{LT2}\cite{Beh}). In other words, $\E_1$, $\E_2$, $\F_1$ and $\F_2$ are global vector bundles, where $\Edot = [\E_1 \to \E_2]$ and $\Fdot = [\F_1 \to \F_2]$ are the standard perfect obstruction theories.

First we recall the notion of compatible perfect obstruction theories (see \cite{LT2}\cite{inc}\cite{functor}\cite{deg}) :
\begin{definition}{\label{defcompatible}}
Suppose $\vp : \bM \to \bN$ is a morphism between separated Deligne-Mumford stacks. Let $\Eup$, $\Fup$ and $\Lup$ be the (dual) perfect obstruction theories for $\bM$, $\bN$ and $\bM / \bN$. We say $\Eup$, $\Fup$ and $\Lup$ are compatible if and only if we have a morphism of distinguished triangles (the bottom row is the triangle of cotangent complexes):
$$
\xymatrix{
\vp^*\Fup \ar[r] \ar[d] & \Eup \ar[r] \ar[d] & \Lup \ar[r]\ar[d] & \vp^* \Fup[1]\ar[d]\\
\vp^*L_{\bN} \ar[r] & L_{\bM} \ar[r] & L_{\bM/ \bN}\ar[r] & \vp^*L_{\bN}[1]
}.
$$
\end{definition}

\begin{remark}
There are different versions of compatibility. One version (\cite{LT2},\cite{inc},\cite{functor}) requires $\Lup$ must come from the relative cotangent complex associated to a local complete intersection morphism of relative Deligne-Mumford type. Here we adapt a broader definition, as used in \cite{deg}.
\end{remark} 
Proposition~\ref{surjective} implies the existence of a relative perfect obstruction theory on $\tM =\moduli_{0,n}(\tX, \tbeta) \to M=\moduli_{0,n}(X, \pi_*\tbeta)$.
\begin{lemma}{\label{triangle}}
If $N_{Z/X}$ is convex, then there exists a distinguished triangle in $\deri(\sO_{\tM})$
$$\Ldot \to \Edot \to \vp^*\Fdot \to \Ldot [1]$$
$$\Ldot=[\sL_1 \to \sL_2] \text{ , where } \sL_i \text{ are locally free sheaves.}$$
\end{lemma}

\begin{proof}
One can always find $\Ldot$ so that $\Ldot \to \Edot \to \vp^*\Fdot \to \Ldot [1]$ is a distinguished triangle. Since $\tM$ has enough locally free sheaves (see~\cite{GP}), we may assume $\Ldot=[\sL_1 \to \sL_2 \to \sL_3] \text{ , where } \sL_i \text{ are locally free sheaves.}$
The associated cohomology long exact sequence is 
$$
\xymatrix{
\sO b_{\tM} \ar[r] \ar@{=}[d]& \vp^*(\sO b_M) \ar@{=}[d]\\ 
h^2(\Edot) \ar[r] & h^2(\vp^*\Fdot) \ar[r] & h^3(\Ldot) \ar[r]&0
.}
$$
By Proposition~\ref{surjective}, we know $h^3(\Ldot)=0$. This implies 
$$\truntwo(\Ldot) \to \Ldot \text{ is quasi-isomorphic, where } \truntwo(\Ldot)=[\sL_1 \to K_2]$$
with the short exact sequence of sheaves
$$0 \to K_2 \to \sL_2 \to \sL_3 \to 0.$$
Because $\sL_2$ and $\sL_3$ are locally free, $K_2$ is locally free as well.
Replace $\Ldot$ by $\truntwo(\Ldot)$ and change the arrow accordingly, this completes the proof.
\end{proof}

\begin{proposition}{\label{compatible}}
$\Ldot$ introduced in Lemma~\ref{triangle} gives rise to compatible perfect obstruction theories on $\vp : \tM \to M$.
\end{proposition}

\begin{proof}
Define 
\xymatrix{
\Lup:= (\Ldot)^{\vee}[-1], & \Eup:= (\Edot)^{\vee}[-1], & \Fup:= (\Fdot)^{\vee}[-1]
}
\\
Note $\Eup$ and $\Fup$ are the perfect obstruction theories used in~\cite{Beh} and~\cite{inc}.
We have a distinguished triangle
$$\vp^*\Fup \to \Eup \to \Lup \to \vp^*\Fup[1].$$
We also have a distinguished triangle of cotangent complexes
$$\vp^*L_M \to L_{\tM} \to L_{\tM/M} \to \vp^*L_M[1].$$
By the axiom of derived categories, we have a morphism of distinguished triangles:
\begin{equation}
\xymatrix{
\vp^*\Fup \ar[r] \ar[d]^{\alpha}& \Eup \ar[r] \ar[d]^{\beta}& \Lup \ar[r]\ar@{-->}[d]^{\gamma}  & \vp^* \Fup[1]\ar[d]^{\alpha[1]}\\
\vp^*L_M \ar[r] & L_{\tM} \ar[r] & L_{\tM/M}\ar[r] & \vp^*L_M[1]
}
\label{trianglemap}
\end{equation}
Take the associated cohomology long exact sequences of Diagram(~\ref{trianglemap}), we obtain 
$$
\xymatrix{
H^{-1}(\vp^*\Fup) \ar[r] \ar[d]^{\text{surjective}} & H^{-1}(\Eup) \ar[r] \ar[d]^{\text{surjective}}& H^{-1}(\Lup) \ar[r] \ar[d]^{h^{-1}(\gamma)}& H^0(\vp^*\Fup) \ar[r] \ar[d]^{\cong}& H^0(\Eup) \ar[r]\ar[d]^{\cong}& H^0(\Lup)\ar[r]\ar[d]^{h^0(r)} &0 \ar[d]\\
H^{-1}(\vp^*L_M) \ar[r] & H^{-1}(L_{\tM}) \ar[r] & H^{-1}(L_{\tM/M}) \ar[r] & H^0(\vp^*L_M) \ar[r] & H^0(L_{\tM}) \ar[r] & H^0(L_{\tM/M}) \ar[r] & 0
.}
$$
By diagram chasing, we know $h^{-1}(\gamma)$ is surjective and $h^0(\gamma)$ is an isomorphism.
\end{proof}
\begin{lemma}
Suppose $\mathcal{A}, \mathcal{B}, \mathcal{C}$ are separated DM-stacks equipped with perfect obstruction theories.\\
If $f: \mathcal{A} \to \mathcal{B}$ and $g: \mathcal{B} \to \mathcal{C}$ both have compatible perfect obstruction theories, then so does the composition map $g\circ f : \mathcal{A} \to \mathcal{C}$.
\end{lemma}
\begin{proof}
This is a consequence of the octahedron axiom.
\end{proof}

\subsection{Comparison of virtual classes}
In this section, we assume $\vp: \bM \to \bN$ is a morphism between separated Deligne-Mumford stacks. All results will be applied to the case $\bM = \moduli_{0,n}(\tX, \pi^!\beta)$ and $\bN=\moduli_{0,n}(X, \beta)$, where $\beta \in H_2(X)$. Note that $\pi_* \pi^! \beta=\beta$.

There are two equivalent approaches to virtual classes \cite{inc}\cite{LT2}\cite{functor}.
On the other hand, when $\vp : \bM \to \bN$ carries compatible perfect obstruction theories, there is also a different construction of the virtual class of $\bM$, as shown in Section 4.1 of~\cite{deg}. The main argument is the associativity of Gysin maps.
We will quote Lemma 4.3 in~\cite{deg} in the following situation:
\begin{proposition}{\label{relative}}
Given a morphism $\vp : \bM \to \bN$ of separated Deligne-Mumford stacks, if $\vp$ carries compatible perfect obstruction theories, then one can construct a class $[\bM, \bN]^\vir$ in $A_*(\bM)$, and we have $[\bM]^\vir=[\bM, \bN]^\vir$ in $A_*(\bM)$.
\end{proposition}

Suppose $\Ldot$, $\Edot$ and $\Fdot$ are compatible perfect obstruction theories on $\bM/\bN$, $\bM$ and $\bN$ respectively. Now we fix notation in the construction of $\bMN$. Define $Ob_{\bM/\bN}:=h^2(\Ldot)$ as the relative obstruction sheaf.
There is an infinitesimal model (denoted by $\Dptwoset$ in~\cite{deg}) over the pair $(\bM, Ob_{\bM/\bN}\oplus \vp^*Ob_{\bN})$. Consider the surjective map
$$ \sL_2\oplus \vp^*\F_2 \to Ob_{\bM/\bN}\oplus \vp^*Ob_{\bN} \to 0,$$
this gives rise to a cone $\C_{\bM} \subset \vect({\sL_2\oplus \vp^*\F_2})$, so that $\C_{\bM}$ is consistent with $\Dptwoset$. The second construction $\bMN$ is defined as the intersection class of $\C_{\bM}$ with the zero section of $\sL_2\oplus \vp^*\F_2$.\\

The "construction of $\bMN$" in the setting of Behrend-Fantechi construction has appeared in Theorem 1 in \cite{functor}. Theorem 1 in~\cite{functor} is only formulated in the case where $\Lup$ is the pull-back of a relative cotangent complex associated to a local complete intersection morphism of relative Deligne-Mumford type. However, the second part in the proof of Theorem 1 in~\cite{functor} doesn't rely on "local complete intersection morphism", therefore the proof can be slightly rearranged to give "the construction of $\bMN$" in the broader defintion of compatibility (Definition~\ref{defcompatible}). Here we briefly describe how this is achieved by the argument in~\cite{functor}.

Let $\sC_{\bN}$ be the (intrinsic) normal cone stack of $\bN$, and let $\sC_{\bM / \bN}$ be the relative normal cone stack of $\bM / \bN$. One can form another normal cone stack $\sC_{\bM / \sC_{\bN}}$, which is a natural subcone stack of $\sC_{\bM / \bN} \times_{\bM} \vp^*\sC_{\bN}$. Therefore $\sC_{\bM / \sC_{\bN}}$ embeds in the vector bundle stack $\sigma : h^2/h^1(\Ldot) \oplus \vp^*(h^2/h^1(\Fdot)) \to \bM$. The new class $\bMN$ is defined as $(\sigma^*)^{-1} ([\sC_{\bM / \sC_{\bN}}])$.

Given a morphism $X \to Y$ of relative Deligne-Mumford type, denote the deformation (to the normal cone) stack by $\mathcal{M}^0_{X/Y} \to \bP ^1$, with the fiber over $ \{ 0 \} \in \bP^1$ isomorphic to the normal cone stack $\sC_{X/Y}$. If $Y= \text{spec}(\rC)$, denote the deformation stack simply by $\mathcal{M}^0_{X}$. In order to show 
$$(\sigma^*)^{-1} ([\sC_{\bM / \sC_{\bN}}]) = [\bM]^\vir \in A_*(\bM),$$
one considers the double deformation stack $\mathcal{M}^0_{\bM \times \bP^1 / \mathcal{M}^0_{\bN}} \to \bP^1 \times \bP^1 $. This provides a rational equivalence $[\sC_{\bM / \sC_{\bN}}] \approx [\sC_{\bM}]$ in $\sC_{\bM \times \bP^1 / \mathcal{M}^0_{\bN}}$.
On the other hand, by Proposition 1 in~\cite{functor}, the abelian hull of $\sC_{\bM \times \bP^1 / \mathcal{M}^0_{\bN}}$ has a natural map to the vector bundle stack $h^1/h^0(c(g))$ on $\bM\times \bP^1$, where $c(g)$ is the mapping cone associated to 
$$\Edot\oplus\vp^*\Fdot \stackrel{g}{\rightarrow} \vp^*\Fdot\boxtimes\sO_{\bP^1}(1) \text{ on } \bM\times \bP^1.$$
Now the rational equivalence can be pushed forward to $h^1/h^0(c(g))$. It is easy to see that the pull back of $h^1/h^0(c(g))$ to $\bM\times \{0\}$ and $\bM\times \{1\}$, correspond to $h^2/h^1(\Ldot) \oplus \vp^*(h^2/h^1(\Fdot))$ and $h^2/h^1(\Edot)$ respectively. Therefore $(\sigma^*)^{-1} ([\sC_{\bM / \sC_{\bN}}]) = [\bM]^\vir$.\\

\begin{remark}
At the beginning of Section 4.1 in~\cite{deg}, it is assumed that $\bM \to \bN$ is representable. One can drop this assumption by taking a presentation of $\bM$: a surjective et\'ale morphism from a scheme $T \to \bM$. There are natural compatible perfect obstruction theories on $T \to \bN$ induced from those on $\bM \to \bN$. Note $T \to \bN$ is representable, so we can apply Lemma 4.3 in~\cite{deg}. On the other hand, the construction of various cones, cycles and rational equivalence in the proof of Lemma 4.3 are canonical, and they descend to the case $\bM \to \bN$. Alternatively, this can also be seen via the construction in~\cite{functor}, as described in the previous paragraph.
\end{remark}

\begin{remark}
The "construction of $\bMN$" is only useful when one has a good understanding of the relative obstruction theory $\Lup$, otherwise it simply transforms a problem into something unknown. In practice, it is usually quoted in the form of Theorem 1 in~\cite{functor}, where $\Lup$ comes from local complete intersection.
\end{remark}

\bigskip

Regarding the construction of $[\bN]^\vir$, take the surjective map \quad
$ \F_2 \to Ob_{\bN} \to 0$.\\ This gives rise to a cone $\C_{\bN} \subset \vect(\F_2)$, so that $\C_{\bN}$ is consistent with the infinitesimal model over $(\bN, Ob_{\bN})$. $[\bN]^\vir$ is defined as the intersection class of $\C_{\bN}$ with the zero section of $\F_2$.
Note that $\C_{\bM}$ is a cone with pure dimension $=\vdim (\bM) + rk(\sL_2) + rk(\F_2)$, and $\C_{\bN}$ is a cone with pure dimension $=\vdim (\bN) + rk(\F_2)$.

\begin{lemma}{\label{cone}}
We have the following diagram (not Cartesian product):
$$\xymatrix{
C_{\bM} \ar[r]\ar[dd] & \vect (\sL_2 \oplus \vp^* \F_2) \ar[d]\\
 & \vect (\vp^* \F_2) \ar[d]\\
C_{\bN} \ar[r] & \vect(\F_2) 
}
$$
\end{lemma}

\begin{proof}
The properties of $\C_{\bM}$ and $\C_{\bN}$ are determined by properties of infinitesimal models, therefore it suffices to prove the corresponding diagram in the infinitesimal models, which is straightforward. An alternative way to see this is via the construction in~\cite{inc} and \cite{functor}.
\end{proof}
On the other hand, $\vectLvF$ can be also regarded as a vector bundle over $\vectvF$. Let 
$$0_{\sL_2}: \vectvF \to \vectLvF$$ be the zero section, then we have the diagram:
$$
\xymatrix{
\C_{\bN} & \ar[l]_(.7){\psi}\C_{\bM} \cap \vectvF \ar[r] \ar[d]& \C_{\bM} \ar[d]\\
& \vectvF \ar[r]^{0_{\sL_2}} & \vectLvF 
}
$$
The right square is a Cartesian product, and $\psi$ is the map induced by $\C_{\bM} \to \C_{\bN}$ in Lemma~\ref{cone}. 
Note $\psi$ is proper as long as $\vp$ is proper. We also have 
$$\zcone \in A_{d + rk(\F_2)}(\C_{\bM} \cap \vectvF) \text{ , where } 0_{\sL_2}^! \text{ is the refined Gysin map.}$$

In the blow-up case $\vp: \moduli_{0,n}(\tX, \pi^!\beta) \to \moduli_{0,n}(X,\beta)$, consider the following diagram:
$$
\xymatrix{
\moduli_{0,n+1}(X,\beta) \ar[r]^>>>>>>{e_{n+1}} \ar[d]^{\pi_{n+1}}& X\\
\moduli_{0,n}(X,\beta)
}
$$
Let $U$ be the complement of $\pi_{n+1} (e_{n+1}^{-1} (Z))$ in $\moduli_{0,n}(X,\beta)$, therefore $U$ is an open substack of $\moduli_{0,n}(X,\beta)$. 
Given $(C,\underline{a},f) \in M$, we have:
$$  (C,\underline{a},f) \in U \Leftrightarrow f(C)\cap Z =\phi.$$
Because $\pi: \tX \to X$ is the blow up of $X$ along $Z$, we deduce:
\begin{lemma}{\label{iso}}
There is an isomorphism $\vp: \vp^{-1}(U) \to U$ with the same (in the sense of quasi-isomorphic) perfect obstruction theory.
\end{lemma}
Because of the above lemma, it motivates us to analyze the following situation:\\
Suppose the proper morphism $\vp : \bM \to \bN$ has compatible perfect obstruction theories with $d=\vdim(\bM)=\vdim(\bN)$. Moreover, we assume that there exists an open substack $U$ in $\bN$, so that $\vp: \vp^{-1}(U) \to U$ is an isomorphism with the same perfect obstruction theories.

\begin{lemma}{\label{openpart}}
Under the setting in the previous paragraph, we have 
$$\left(\psi_*\zcone\right)|_U = [\C_{\bN}|_U] \text{ in } A_{d + rk(\F_2)}(\C_{\bN}|_U), \text{ where } d=\vdim{\bM}=\vdim{\bN}$$
\end{lemma}

\begin{proof}
Because $\vp^{-1}U \cong U$ and the flat pull back $(\bullet)|_U$ commutes with other operators, we have
$$\left(\psi_*\zcone\right)|_U = \psi_* 0_{\sL_2}^!(\C_{\bM}|_U).$$
Moreover, according to $\Edot|_{\vp^{-1}(U)}\cong \Fdot|_U$, we know $Ob_{\bM/\bN}|_U$ vanishes. Therefore the infinitesimal models on $\vp^{-1}(U) \cong U$ are the same, and we have the Cartesian diagram:
$$
\xymatrix{
\C_{\bM}|_U \ar[r] \ar[d] & \vect(\sL_2|_U \oplus \F_2|_U)\ar[d]\\
\C_{\bN}|_U\ar[r] & \vect(\F_2|_U)
}
$$
In other words, $\vect(\sL_2|_U \oplus \F_2|_U)$ is a vector bundle over $\vect(\F_2|_U)$, and $\C_{\bM}|_U$ is the flat pull back of $\C_{\bN}|_U$.
Therefore we have $\psi_* 0_{\sL_2}^!(\C_{\bM}|_U) = [\C_{\bN}|_U]$ in $A_{d + rk(\F_2)}(\C_{\bN}|_U)$.
\end{proof}

Suppose $\C_{\bN}$ has irreducible components $\C_i$, $i = 1, \cdots, k$. Let $\supp(\bullet)$ be the support of a cone.
In the rest of this section, we will assume the open substack $U \subset \bN$ satisfies the following technical assumption:
\begin{equation}\tag{Assumption $\ast$}\label{cond*}
\supp(\C_i) \cap U \text{ is non-empty in } \bN \text{ for } i = 1, \cdots, k.
\end{equation}
Because $\C_{\bN}$ is a cone with pure dimension $d + rk(\F_2)$, (~\ref{cond*}) implies
$$\dim (\C_{\bN} - (\C_{\bN}|_U )) < d+\rk(\F_2) .$$
With this technical assumption, it is easy to prove that $\vp : \bM \to \bN$ is virtually birational.

\begin{corollary}{\label{basiclemma}}
If the open substack $U \subset \bN$ satisfies~\ref{cond*}, then we have
$$\psi_*\zcone = [\C_{\bN}] \text{ in } A_{d + rk(\F_2)}(\C_{\bN}).$$
\end{corollary}

\begin{proof}
We have
$$\dim (\C_{\bN} - (\C_{\bN}|_U )) < d+\rk(\F_2) \Longrightarrow A_{d + rk(\F_2)}(\C_{\bN}-\C_{\bN}|_U)=0.$$
Combined with the right exact sequence:
$$A_{d + rk(\F_2)}(\C_{\bN}-\C_{\bN}|_U) \to A_{d + rk(\F_2)}(\C_{\bN}) \to A_{d + rk(\F_2)}(\C_{\bN}|_U) \to 0,$$
we know $A_{d + rk(\F_2)}(\C_{\bN}) \cong A_{d + rk(\F_2)}(\C_{\bN}|_U)$. By Lemma~\ref{openpart}, we obtain $$\psi_*\zcone = [\C_{\bN}] \text{ in } A_{d + rk(\F_2)}(\C_{\bN}).$$
\end{proof}
Now we summarize all results in this section to deduce the following:

\begin{proposition}{\label{special}}
Suppose the proper morphism $\vp : \bM \to \bN$ has compatible perfect obstruction theories with $d=\vdim(\bM)=\vdim(\bN)$. We also assume that there exists an open substack $U$ in $\bN$, so that $\vp: \vp^{-1}(U) \to U$ is an isomorphism with the same perfect obstruction theories.\\\\
If $U \cap \supp (\text{ each irreducible component of }\C_{\bN}) \text{ is non-empty in } \bN$,\\
then we have $\vp_* [\bM]^\vir = [\bN]^\vir$ in the Chow group $A_d(\bN)$.
\end{proposition}

\begin{proof}
Recall the diagram with the right square as Cartesian product:
$$
\xymatrix{
\C_{\bN} \ar[d] & \vectvF \cap \C_{\bM} \ar[l]_(.6){\psi} \ar[r] \ar[d]& \C_{\bM}\ar[d]\\
\vect(\F_2) & \vectvF \ar[l]_{\vp} \ar[r]^(.4){0_{\sL_2}} & \vectLvF
}$$
If we regard $[\C_{\bM}]$ as a class in $A_*(\vectLvF)$, and $[\C_{\bN}]$ as a class in $A_*(\vectF)$, then Corollary~\ref{basiclemma} implies $$\vp_*\zcone = [\C_{\bN}] \text{ in }A_{d + rk(\F_2)}(\vectF).$$
From another diagram with the left square as Cartesian product:
$$
\xymatrix{
\bM \ar[rr]^(.4){0_{(\vp^*\F_2)}} \ar[d]_{\vp} && \vectvF \ar[r]^<<<<<<{0_{\sL_2}} \ar[d]^{\vp} & \vectLvF\\
\bN \ar[rr]^(.4){0_{\F_2}} & & \vect(\F_2)
}
$$
\begin{align*}
\vp_*[\bM]^{\vir} &= \vp_* \circ 0^!_{(\vp^*\F_2)} \circ \zcone = 0^!_{\F_2} \circ \vp_* \circ \zcone\\
&= 0^!_{\F_2}[\C_{\bN}] = [\bN]^{\vir}
\end{align*}
Here we use $0_{(\vp^*\F_2)}^!= 0_{\F_2}^! : A_* (\vect (\vp^*\F_2)) \to A_*(\bM)$.
\end{proof}

\begin{corollary}{\label{stupid}}
Suppose the proper morphism $\vp : \bM \to \bN$ has compatible perfect obstruction theories with $d=\vdim(\bM)=\vdim(\bN)$. Suppose there exists an open substack $U$ in $\bN$, so that $\vp(\bM)\cap U = \emptyset$.\\\\
If $U \cap \supp (\text{ each irreducible component of }\C_{\bN}) \text{ is non-empty in } \bN$,\\
then we have $\vp_* [\bM]^\vir = 0$ in the Chow group $A_d(\bN)$.
\end{corollary}
\begin{proof}
Apply the previous proposition to $\vp \coprod \text{Id} : \bM \coprod \bN \to \bN$.
\end{proof}
\subsection{Transversal intersection of two manifolds}
Suppose $X$ and $Y$ are two arbitrary closed submanifolds of a compact homogeneous space $\mathcal{P}$, and $Z$ is the transversal intersection of $X$ and $Y$. Suppose the group variety $G$ acts on $\mathcal{P}$ transitively.
\begin{lemma}
The normal bundle $N_{Z/X}$ is generated by global sections, and therefore is convex.
\end{lemma}
\begin{proof}
The tangent bundle $T\mathcal{P}$ is generated by global sections, and $N_{Y/\mathcal{P}}$ is a quotient bundle of $T\mathcal{P}$.
This implies $N_{Y/\mathcal{P}}$ is generated by global sections as well. Note $N_{Z/X}$ is the pull back of $N_{Y/\mathcal{P}}$ to $Z$.
\end{proof}
Consider $\pi: \tX \to X$, the blow up of $X$ along the submanifold $Z$. May assume $\codim(Y, \mathcal{P}) \geq 2.$ The first attempt is to apply Proposition~\ref{special}, but the technical assumption 
$$\supp(\text{ each irreducible component of }\C_{\bN}) \cap U \text{ is non-empty in } \bN, \text{ where } \bN=\moduli_{0,n}(X,\beta)$$
may not be satisfied. We will choose an element $\sigma\in G$, and show the technical assumption is satisfied when Z is perturbed to $X\cap Y^{\sigma}$.

\begin{lemma}
Given a holomorphic map from a compact curve $f: C \to \mathcal{P}$, define 
$$B_{(C,f)}:=\{\sigma\in G \quad | \quad f(C)\cap Y^{\sigma} \neq \emptyset   \}.$$
Then $B_{(C,f)}$ is closed in $G$, and $\dim (B_{(C,f)}) < \dim G$.
\end{lemma}
\begin{proof}
Consider $G \stackrel{p_1}{\longleftarrow} G\times Y \stackrel{\Phi}{\longrightarrow} \mathcal{P}$, where $p_1$ is the projection and $\Phi$ is the group action. Note $B_{(C,f)}= p_1(\Phi^{-1}(f(C)))$ is closed in $G$ because $p_1$ is proper. Because $G$ acts on $\mathcal{P}$ transitively, $\Phi$ is a smooth morphism. Therefore, 
$$\dim B_{(C,f)} \leq \dim Y + \dim G-\dim \mathcal{P}+ \dim f(C) \leq \dim G - \codim (Y, \mathcal{P})+1 \leq \dim G-1.$$
\end{proof}

\begin{lemma}{\label{dimtrans}}
Define $W:=\{\sigma \in G \quad |\quad Y^{\sigma} \text{ is not transversal to } X \}$. Then $W$ is closed in $G$,\\ with $\dim W < \dim G$.
\end{lemma}
\begin{proof}
Note the identity element $\text{ Id} \in G$ does not belong to $W$.
\end{proof}

Suppose $\Fdot$ is the perfect obstruction theory on $\bN=\Monxb$, and the virtual normal cone $\C_{\bN}$ has irreducible components $\C_i$, for $i=1, \dots, k$. For each $i$, we pick a point $(C_i, \underline{a}_i, f_i)\in \text{supp}(\C_i)$. By the previous two lemmas, we know
$$ W\cup (\bigcup^k_{i=1} B_{(C_i, f_i)}) \text{ is closed in } G \text{ with codimension } \geq 1.$$
Therefore we can take an affine smooth locally closed curve $S \hookrightarrow G$ such that:
\begin{enumerate}
\item $\text{Id} \in S$,
\item $(S-\text{Id}) \cap (\bigcup^k_{i=1} B_{(C_i, f_i)}) = \emptyset$,
\item $Y^{\sigma}$ is transversal to $X$, $\forall \sigma\in S$.
\end{enumerate}

Choose an element $\sigma \in S$, then $Z= X\cap Y$ is deformation equivalent to $Z_{\sigma}:=X \cap Y^{\sigma}$. Note the normal bundle $N_{Z_{\sigma} /X }$ is still generated by global sections.
The technical assumption of Proposition~\ref{special} is satisfied for $Bl_{Z_{\sigma}} X \to X$.

\begin{lemma}
$$\supp(\text{ each irreducible component of }\C_{\bN}) \cap U \text{ is non-empty in } \bN=\Monxb,$$
where $U$ is a collection of stable maps supported away from $Z_{\sigma}$.
\end{lemma}
\begin{proof}
The choice of the curve $S \hookrightarrow G$ asserts $f_i(C_i) \cap Y^{\sigma} =\emptyset$, for $i=1, \dots, k$. Hence $(C_i, \underline{a}_i, f_i) \in U$.
\end{proof}

Proposition~\ref{special} and deformation invariance of virtual classes implies:
\begin{theorem}{\label{trans}}
Suppose $Z$ is the transversal intersection of two manifolds $X$ and $Y$ in a compact homogeneous space $\mathcal{P}$.
Then we have $\vp_*[\overline{\mathcal{M}_{0,n}}(Bl_Z X, \pi^!\beta)]^{\vir} = [\Monxb]^{\vir}$ in the Chow group.
\end{theorem}

The theorem can be generalized to the case when $Z$ is the intersection of $X$ with multiple manifolds.
More precisely, suppose $Y_i$, $i=1, \dots, m$, are submanifolds of a homogeneous space $P$. We assume $Y_{k+1}$ is transversal to $X\cap (\bigcap_{i=1}^k Y_i)$, for $k=0, \dots , m-1$.
\begin{corollary}
Let $Z$ be $X\cap (\bigcap_{i=1}^m Y_i)$, then we have 
$$\vp_*[\overline{\mathcal{M}_{0,n}}(Bl_Z X, \pi^!\beta)]^{\vir} = [\Monxb]^{\vir}.$$
\end{corollary}
\begin{proof}
Define $G':= G^m$ and $\mathcal{P}':= \mathcal{P}^m$. Then $\mathcal{P}'$ is a homogeneous space with respect to the group variety $G'$. Let $\Delta : X \to \mathcal{P}'$ be the diagonal map. $X$ is transversal to the manifold $\prod_{i=1}^m Y_i$ in the ambient space $\mathcal{P}'$. Apply the previous theorem to the case $X \to \mathcal{P}'$, $Y:=\prod_{i=1}^m Y_i$, and $Z= X \cap Y$ in $\mathcal{P}'$.
\end{proof}
A similar argument also implies:
\begin{corollary}
Suppose $X$ is a projective manifold, and $Z$ is a collection of points in $X$.
Then $$\vp_*[\overline{\mathcal{M}_{0,n}}(Bl_Z X, \pi^!\beta)]^{\vir} = [\Monxb]^{\vir}.$$
\end{corollary}
\begin{proof}
This is because $N_{Z/X}$ is convex, and $Z$ can always be moved.
\end{proof}

\section{Virtual Birationality after degeneration}
In the previous subsection, the submanifold is deformed so that the technical assumption in Proposition~\ref{special} is satisfied. In general, if $N_{Z/X}$ has a non-zero section, it doesn't imply $Z$ can be moved. Degeneration formula reduces the problem to a ruled variety, where $Z$ can be moved if $N_{Z/X}$ has a section.

Degeneration formula has been clearly presented in~\cite{deg} \cite{relloc}, and \cite{LiuBlowup}. The purpose of the first subsection is to fix notation.

\subsection{Degeneration formula from blow-ups}
$(X, D)$ is called a relative pair if $D$ is a smooth divisor of the manifold $X$.
The $\bP ^1$-bundle $Y:=\bP_{D}(N_{D/X}\oplus \sO_D)$ has two divisors: the zero divisor (with normal bundle $N_{D/X}^\vee$) and the infinity divisor (with normal bundle $N_{D/X}$). $Y_l$ is defined as the union of $l$-copies of $Y$, by gluing the infinity divisor of the $i$-th component to the zero divisor of the $(i+1)$-th component. Let $D_i$ be the zero divisor of the $(i+1)$-th component. Note $\text{Sing}(Y_l)=\cup_{i=1}^{l-1} D_i$. Define $X_l:= X \cup_{D_0} Y_l$.\\
Now we recall Definition 4.6 in~\cite{deg}:
\begin{definition}
An admissible weighted graph $\Gamma$ for a relative pair $(X, D)$ is a graph without edges together with the following data:
\begin{enumerate}
\item an ordered collection of legs, an ordered collection of weighted roots, and two weight functions on the vertex set $g: V(\Gamma) \to \mathbb{Z}_{\geq 0}$ and $b: V(\Gamma) \to H_2(X)$.
\item $\Gamma$ is relatively connected in the sense that either $|V(\Gamma)|=1$ or each vertex in $V(\Gamma)$ has at least one root attached to it.
\end{enumerate}
\end{definition}

The weight functions $g$ and $b$ in the previous defintion are used to record the genus and the homology class in $X$ for each connected component of a domain curve.
Denote the moduli of relative stable maps to $(X, D)$ with topological data $\Gamma$ by $\overline{\mathcal{M}}(X, D, \Gamma)$. A $\mathbb{C}$-point in $\overline{\mathcal{M}}(X, D, \Gamma)$ is a holomorphic map $C \stackrel{f}{\rightarrow} X_l \to X$ satisfying \textbf{stability} and \textbf{predeformability} together with topological constraints $\Gamma$.
The domain curve is disconnected if and only if $|V(\Gamma)|>1$. For convenience, define
$$b(\Gamma) := \sum_{v \in V(\Gamma)} b(v) \text{   and   } g(\Gamma):= 1- |V(\Gamma)| + \sum_{v \in V(\Gamma)} g(v).$$
Let $\mathcal{T}$ be the Artin stack parametrizing the possible target of relative stable maps to $(X, D)$. The perfect obstruction theory on $\moduli (X, D, \Gamma)$ is induced from the relative perfect obstruction theory on 
$$\moduli (X, D, \Gamma) \to \mathcal{T} \times \mathfrak{M}_{g(\Gamma), k}, \text{ where } k=\text{ number of legs } + \text{ number of roots}.$$

Given an arbitrary manifold $X$ with a submanifold $Z$, deformation to the normal cone is obtained from the blow-up of a trivial family: $$ W:= \text{Bl}_{Z\times \{ 0 \} } X \times \mathbb{C} \to \mathbb{C}.$$
Note $W_t \cong X$ for $t \neq 0$ and $W_0= \tX \bigcup _ {\bP(N_{Z/X})} \bP(N_{Z/X} \oplus \sO_Z)$, where $\pi : \tX= \text{Bl}_Z X \to X$. Denote $\bP(N_{Z/X})$ by $D$, and $p_2: \bP(N_{Z/X} \oplus \sO_Z) \to Z$.

\begin{theorem}[Degeneration formula from blow-up, see~\cite{deg} and \cite{LiuBlowup}]

$$ [\moduli (\mathcal{W}_0, g, n , \beta)]^\vir = \sum_{\eta \in \overline{\Omega}_{(g, n, \beta)}} \frac{\textbf{m}(\eta)}{|\text{Eq}(\eta)|} \Phi_{\eta_*} \Delta ^! ([\moduli (\tX, D, \Gamma_1)]^\vir \times [\moduli (\bP(N_{Z/X} \oplus \sO_Z), D, \Gamma_2)]^\vir), \text{ where } \beta \in H_2(X).$$
\end{theorem}

The set $\overline{\Omega}_{(g, n, \beta)}$ is an equivalence set $\Omega_{(g, n, \beta)} / \sim_{\text{equ}}$. The set $\Omega_{(g, n, \beta)}$ is a collection of admissible triples $\eta = (\Gamma_1, \Gamma_2, I)$ satifying:
\begin{enumerate}
\item $\Gamma_1$ and $\Gamma_2$ are admissible weighted graphs for $(X, D)$ and $(\bP(N_{Z/X} \oplus \sO_Z), D)$ respectively.
\item $\Gamma_1$ and $\Gamma_2$ are required to have identical number of roots, say $r$ roots. The weight of $i$-th root in $\Gamma_1$ and $\Gamma_2$ must be identical, for $i=1, \cdots, r$. 
\item If one glues all corresponding roots of $\Gamma_1$ and $\Gamma_2$, then the new graph must be connected.
\item $n= \# \text{legs}(\Gamma_1) + \# \text{legs}(\Gamma_2).$
\item $I$ is a rule concerning the ordering of the union of legs in $\Gamma_1$ and $\Gamma_2$.
\item (Genus constraint) $g(\eta) := g(\Gamma_1) + g(\Gamma_2) + r -1$ must equal $g$.
\item (Homology constraint) $\pi_*(b(\Gamma_1)) + p_{2_*}(b(\Gamma_2)) = \beta$ and some other restrictions, see Section 3 in~\cite{LiuBlowup}.
\end{enumerate}

Given a permutation $\sigma \in S_r$, $\eta^\sigma$ is defined by reordering $r$ roots. Define $\eta_1 \sim \eta_2$ if and only if $\eta_1^\sigma = \eta_2$ for some $\sigma$. $\overline{\Omega}_{(g, n, \beta)}$ is the equivalence class of this relation.
Define
$$\text{Eq}(\eta):= \# \{\sigma \in S_r | \eta^\sigma = \eta   \} \text{ and } \textbf{m}(\eta):= \text{ the product of the weights of the roots in } \Gamma_1 .$$

$$
\xymatrix{
\moduli (\tX, D, \Gamma_1) \times_{D^r} \moduli (\bP(N_{Z/X} \oplus \sO_Z), D, \Gamma_2) \ar[d] \ar[r] & \moduli (\tX, D, \Gamma_1) \times \moduli (\bP(N_{Z/X} \oplus \sO_Z), D, \Gamma_2) \ar[d]\\
D^r \ar[r]^{\Delta \text{ diagonal}} & D^r \times D^r
}
$$
$\Phi_{\eta} : \moduli (\tX, D, \Gamma_1) \times_{D^r} \moduli (\bP(N_{Z/X} \oplus \sO_Z), D, \Gamma_2) \to \moduli (\mathcal{W}_0, g, n , \beta)$ is to glue two relative stable morphisms.

\bigskip

One can also apply the deformation to the normal cone to $D \subset \tX$:
$$ \tilde{W}:= \text{Bl}_{D\times \{ 0 \} } \tX \times \mathbb{C} \to \mathbb{C}.$$
Note $\tilde{W}_t \cong \tX$ for $t \neq 0$ and $\tilde{W_0}= \tX \bigcup _ {\bP(N_{Z/X})} \overline{\bP(N_{Z/X} \oplus \sO_Z)}$, where $\overline{\bP(N_{Z/X} \oplus \sO_Z)}$ is the blow up of $\bP(N_{Z/X} \oplus \sO_Z)$ along $Z$. This space can also be viewed as a $\bP^1$-bundle over $\bP(N_{Z/X})$:
$$\tilde{p_2}: \overline{\bP(N_{Z/X} \oplus \sO_Z)} 
= \bP_D(\sO_{N_{Z/X}}(1) \oplus \sO) \to \bP(N_{Z/X})=D.$$
Our goal is to compare the virtual classes $[\moduli (\mathcal{W}_0, 0, n , \beta)]^\vir$ and
$[\moduli (\tilde{\mathcal{W}_0}, 0, n , \pi^! \beta)]^\vir$. By the degeneration formula, the main issue is to realize all contributions from $(\bP(N_{Z/X} \oplus \sO_Z), D)$ and $(\overline{\bP(N_{Z/X} \oplus \sO_Z)}, D)$.
\subsection{Relative case}
Let $\pi_Y : \tY \to Y:=\bP(N_{Z/X} \oplus \sO_Z)$ be the blow up along $Z$. Given an adimissible graph $\tGamma$ for $(\tY, D)$, define the adimissible graph $\pi_{Y_*}(\tGamma)$ for $(Y, D)$ by the following:
\begin{enumerate}
\item All information of $\pi_{Y_*}(\tGamma)$ is identical to that of $\tGamma$ except the weight function $b$.
\item We have a commutative diagram:
$$
\xymatrix{
V(\tGamma) \ar[d] \ar[r]^{b} & H_2(\tY) \ar[d]\\
V(\pi_{Y_*}(\tGamma)) \ar[r]^{b} & H_2(Y)
}
$$
\end{enumerate}

\begin{lemma}{\label{homology}}
Suppose $N_{Z/X}$ is convex, and the genus weight function of $\tGamma$ is a zero function. Assume $\moduli (Y, D, \pi_{Y_*}(\tGamma))$ makes sense.\\
Then $\moduli (\tY, D, \tGamma) \to \moduli (Y, D, \pi_{Y_*}(\tGamma))$ has compatible perfect obstruction theories.
\end{lemma}

\begin{proof}
Let $\mathcal{T}$ be the Artin stack parametrizing the possible target of relative stable maps to $(Y, D)$ and $(\tY, D)$. The perfect obstruction theory on $\moduli (Y, D, \pi_{Y_*}(\tGamma))$ is induced from a relative perfect obstruction theory on 
$$\moduli (Y, D, \pi_{Y_*}(\tGamma)) \to \mathcal{T} \times \mathfrak{M}_{g(\tGamma), k}, \text{ where } k=\text{ number of legs } + \text{ number of roots}.$$
Since $\mathcal{T} \times \mathfrak{M}_{g(\tGamma), k}$ is a smooth Artin stack, we have a right exact sequence:
$$\text{RelOb}(f) \to \text{Ob}(C,f) \to 0,$$ 
where $\text{Ob}(C,f)$ refers to the obstruction space of $\moduli (Y, D, \pi_{Y_*}(\tGamma))$ at the point $ [C \stackrel{f}{\rightarrow} Y_l \to Y]$.
\begin{align*}
\phi : & \moduli (\tY, D, \tGamma ) \to \moduli (Y, D, \pi_{Y_*}(\tGamma))\\
& [\tilde{C} \stackrel{\tilde{f}}{\rightarrow} \tY_l \to \tY] \mapsto [C \stackrel{f}{\rightarrow} Y_l \to Y]
\end{align*}
We have a diagram of right exact sequence:
$$
\xymatrix{
\text{RelOb}(\tilde{f}) \ar[r] \ar[d] & \text{Ob}(\tilde{C},\tilde{f}) \ar[r] \ar[d] & 0\\
\text{RelOb}(f) \ar[r] & \text{Ob}(C,f) \ar[r] & 0
}
$$
\textbf{\underline {Step 1}} \qquad $\text{RelOb}(\tilde{f}) \to \text{RelOb}(f)$ is surjective.\\

There is a natural diagram of exact sequences:
$$
\xymatrix{
0 \ar[r] & H^1 \big( \tilde{C}, \tilde{f}^* T_{\tY_l}(-\text{log} D_{\infty}) \big) \ar[d] \ar[r] & \text{RelOb}(\tilde{f}) \ar[r] \ar[d] & H^0 \big( \tilde{C}, \tilde{f}^{-1} \mathcal{E}xt^1(\Omega_{\tY_l}(\text{log} D_{\infty}), \sO_{\tY_l})    \big) \ar[r] \ar[d] & 0\\
0 \ar[r] & H^1 \big( C, f^* T_{Y_l}(-\text{log} D_{\infty}) \big) \ar[r] & \text{RelOb}(f) \ar[r] & H^0 \big( C, f^{-1} \mathcal{E}xt^1(\Omega_{Y_l}(\text{log} D_{\infty}), \sO_{Y_l})    \big) \ar[r] & 0
}
$$
$\mathcal{E}xt^1(\Omega_{\tY_l}(\text{log} D_{\infty}), \sO_{\tY_l})$ and 
$\mathcal{E}xt^1(\Omega_{Y_l}(\text{log} D_{\infty}), \sO_{Y_l})$ are both supported on $\text{Sing}(\tY_l)=\text{Sing}(Y_l)=\cup_{i=0}^{l-1}D_i$, and these two sheaves are canonically isomorphic to each other. Therefore, it remains to show the first vertical arrow is surjective.
We also have another exact sequence
$$0 \to T_{\tY_l}(-\text{log} D_{\infty}) \to T_{Y_l}(-\text{log} D_{\infty}) \to Q_{N_{Z/X}} \to 0,$$
where $Q_{N_{Z/X}}$ is the universal quotient bundle on $\bP_Z(N_{Z/X})$. Now the proof proceeds as the second part of the proof in Lemma~\ref{surjective}. This concludes Step 1.\\
\textbf{\underline {Step 2}}\\
$\text{RelOb}(\tilde{f}) \to \text{RelOb}(f)$ is surjective $\Longrightarrow$ $\text{Ob}(\tilde{C}, \tilde{f}) \to \text{Ob}(C, f)$ is surjective.\\
By Lemma~\ref{triangle} and Proposition~\ref{compatible}, there exists a relative perfect obstruction theory on $\moduli (\tY, D, \tGamma) \to \moduli (Y, D, \pi_{Y_*}(\tGamma))$.
Moreover, it is compatible with two existing obstruction theories on the two moduli spaces.
\end{proof}

There is an induced map on adimissible triples:
$\Psi : \Omega_{(0, n, \pi^!\beta)} \to \Omega_{(0, n, \beta)}$, where $\Psi(\Gamma_1, \Gamma_2, I)= (\Gamma_1, \pi_{Y_*}(\Gamma_2), I)$.

\begin{lemma}
Suppose $\Psi(\Gamma_1, \Gamma_2, I)= (\Gamma_1, \Gamma_3, I)$, then we have $b(\Gamma_2)=\pi_Y^!(b(\Gamma_3)) \in H_2(\tY).$
\end{lemma}
\begin{proof}
Since $\pi_{Y_*}(b(\Gamma_2))= b(\Gamma_3)$, it suffices to prove  $b(\Gamma_2)\bullet D_{\infty}=0$ in $\tY$, where $\imath_{\infty}: D_{\infty}=\bP(N_{Z/X}) \hookrightarrow \tY$ has normal bundle $\sO_{N_{Z/X}}(-1)$. Let $\imath_0: D_0=\bP(N_{Z/X}) \hookrightarrow \tY$ be the divisor which has normal bundle $\sO_{N_{Z/X}}(1)$.

We have $b(\Gamma_2)=\imath_{0_*}(\tilde{p_2}_*b(\Gamma_2))+ f$, where $f$ is a multiple of the fiber class of $\tilde{p_2}$. It remains to show $f=0$.
$$\tY \stackrel{\tilde{p_2}}{\rightarrow} D=D_0 \stackrel{\imath_0}{\rightarrow} \tY.$$
$(\Gamma_1, \Gamma_2, I) \in \Omega_{(0, n, \pi^!\beta)}$ implies:
$$ \pi^!\beta= b(\Gamma_1) + \tilde{p_2}_*b(\Gamma_2),$$
$$b(\Gamma_1)\bullet D \text{ in } \tX = b(\Gamma_2)\bullet D_0 \text{ in } \tY.$$
Therefore we have
\begin{align*}
0= \pi^!\beta \bullet D 
& = \big( b(\Gamma_1)\bullet D \text{ in } \tX \big) + \big( \tilde{p_2}_*b(\Gamma_2)\bullet D \text{ in } \tX \big)\\ 
& = \big( b(\Gamma_2)\bullet D_0 \text{ in } \tY \big) - \big( \imath_{0_*}\tilde{p_2}_*b(\Gamma_2)\bullet D_0 \text{ in } \tY \big)\\
& = f\bullet D_0 \text{ in } \tY 
\end{align*}
\end{proof}
Given $(\Gamma_1,\Gamma, I) \in \Omega_{(0, n, \beta)}$, define 
$$\Psi^{-1}(\Gamma)=\{ \Gamma_2 \text{ is a admissible weighted graph for } (\tY, D) |  (\Gamma_1, \Gamma_2, I)\in \Omega_{(0, n, \pi^!\beta)} \text{ such that } \Psi(\Gamma_1, \Gamma_2, I)= (\Gamma_1,\Gamma, I) \}.$$
It is straightforward to check that $\Psi^{-1}(\Gamma)$ depends on $(0,n, \beta)$, but is independent of $\Gamma_1$ and $I$.
Given $\tilde{\Gamma} \in \Psi^{-1}(\Gamma)$, we have a natural map $\moduli (\tY, D, \tilde{\Gamma}) \to \moduli (Y, D, \Gamma)$. Note that $\vdim \moduli (\tY, D, \tilde{\Gamma}) = \vdim \moduli (Y, D, \Gamma)$, and the weight functions $g$ of $\tilde{\Gamma}$ and $\Gamma$ are both zero functions.

On the other hand, there is a canonical pre-image $\pi_Y^!(\Gamma)\in \Psi^{-1}(\Gamma)$, which is characterized by:

\begin{enumerate}
\item All information of $\pi_Y^!(\Gamma)$ is identical to that of $\Gamma$ except the weight function $b$.
\item We have a commutative diagram:
$$
\xymatrix{
V(\pi_Y^!(\Gamma)) \ar[d] \ar[r]^{b} & H_2(Y) \ar[d]^{\pi_Y^!}\\
V(\Gamma) \ar[r]^{b} & H_2(\tY)
}
$$
\end{enumerate}

\bigskip

We will consider two classes of submanifolds. The first one is:

\begin{definition}
A connected submanifold $Z \subset X$ is of type I, if the following two conditions are satisfied:
\begin{enumerate}
\item $N_{Z/X}$ is a convex bundle over Z,
\item There is a subbundle $\calF$ in $N_{Z/X}$ with rank $\rk (\calF)\geq 2$, and $\calF$ is generated by global sections.
\end{enumerate}
\end{definition}

For example, $Z \subset X$ is of type I if $N_{Z/X}$ is generated by global sections.

\begin{lemma}{\label{relbirational}}
If $Z \subset X$ is of type I,
then we have 
\begin{enumerate}
\item
$\phi_*[\moduli (\tY, D, \pi_Y^!(\Gamma))]^{\vir} = [\moduli (Y, D, \Gamma)]^{\vir}.$
\item 
$\phi_*[\moduli (\tY, D, \tGamma)]^{\vir} =0$ when $\pi_Y^!(\Gamma) \neq \tGamma \in \Psi^{-1}(\Gamma)$.
\end{enumerate}
\end{lemma}
\begin{proof}
For the first statement, the submanifold $Z$ will be moved so that the technical assumption in Proposition~\ref{special} is satisfied:
$$\supp(\text{ each irreducible component of }\C_{\bN}) \cap U \text{ is non-empty in } \bN,$$
where $\bN=\moduli (Y, D, \Gamma)$ and $U$ is a collection of relative stable maps supported away from the submanifold $Z$.

For each irreducible component of $\C_{\bN}$, we pick a point in the support of the cone
$$ C_i \stackrel{f_i}{\rightarrow} Y_{l_i} \to Y, \text{ for } i= 1, \cdots, k.$$
Since the subbundle $\calF$ is generated by global sections, we have $\oplus_s \sO_Z \to \calF \to 0$.
$$\rC^s \stackrel{p \text{ projection}}{\longleftarrow} Z\times \rC^s = \vect (\oplus_s \sO_Z) \stackrel{\theta \text{ smooth}}{\longrightarrow} \vect (\calF) \stackrel{closed}{\hookrightarrow} \vect (N_{Z/X}) \stackrel{open}{\hookrightarrow} Y.$$
Because $p$ is proper, $p\big( \theta^{-1} (\vect(\calF)\cap f_i(C_i))   \big)$ is closed with dimension $\leq 1+s - \rk(F) \leq s-1$. There exists a point $q\in \rC^s$ such that 
$q \notin p\big( \theta^{-1} (\vect(\calF)\cap f_i(C_i))   \big)$, for all $i$.\\
$q$ induces a section of $N_{Z/X} \to Z$, say $q(Z) \subset \vect (N_{Z/X})$. We have $q(Z) \cap f_i(C_i)= \emptyset$. Move the submanifold $Z$ to $q(Z)$, and notice that the technical assumption is satisfied for the case $\text{Bl}_{q(Z)} Y \to Y$. By Proposition~\ref{special}, we obtain $\phi_*[\moduli (\tY, D, \pi_Y^!(\Gamma))]^{\vir} = [\moduli (Y, D, \Gamma)]^{\vir}.$

For the second statement, the argument is the same, but one applies Corollary~\ref{stupid} instead.
\end{proof}
\begin{proposition}{\label{typeI}}
Suppose $Z \subset X$ is of type I.
Then we have $\phi_*[\moduli (\tilde{\mathcal{W}_0}, 0, n , \pi^! \beta)]^\vir 
= [\moduli (\mathcal{W}_0, 0, n , \beta)]^\vir$.
\end{proposition}
\begin{proof}
By Lemma~\ref{relbirational} and Degeneration formula from blow-up, it remains to check
$$\text{Eq}(\Gamma_1, \pi_Y^!(\Gamma), I)= \text{Eq}(\Gamma_1, \Gamma, I) \text{ and } \textbf{m}(\Gamma_1, \pi_Y^!(\Gamma), I)= \textbf{m}(\Gamma_1, \Gamma, I), \forall (\Gamma_1, \Gamma, I) \in \overline{\Omega}_{(0, n, \beta)},$$
which is straightforward.
\end{proof}
\begin{definition}
A connected submanifold $Z \subset X$ is of type II if every holomorphic map $f : \bP^1 \to Z$ must be a constant map.
\end{definition}

Manifolds of type II is a very limited class of manifolds. Examples include 
\begin{enumerate}
\item higher genus curves, abelian varieties.
\item a fibration with fibers and the base of type II (e.g. product),
\item a submanifold of a manifold of type II.
\end{enumerate}

\begin{proposition}{\label{typeII}}
Suppose $Z \subset X$ is of type II.
Then we have $\phi_*[\moduli (\tilde{\mathcal{W}_0}, 0, n , \pi^! \beta)]^\vir 
= [\moduli (\mathcal{W}_0, 0, n , \beta)]^\vir$.
\end{proposition}
\begin{proof}
Due to the property of $Z$, any vector bundle over $Z$ is automatically convex.
It suffices to prove Lemma~\ref{relbirational} for type II.
First one observes that there is a natural fibration
$$ \moduli (Y,D,\Gamma) \to Z \text{ with nonsingular fibers }\cong \moduli (\bP^m, \bP^{m-1}, \Gamma) , \text{ where } m= \rk(N_{Z/X}).$$
In particular, $\moduli (Y,D,\Gamma)$ is a smooth DM-stack. Therefore the technical assumption of Proposition~\ref{special} is equivalent to saying :
\begin{center}
any point in $\moduli (Y,D,\Gamma)$ can be moved so that the corresponding curve is supported away from $Z$.
\end{center}
The point will be moved along the fiber $\moduli (\bP^m, \bP^{m-1}, \Gamma)$, so we may assume $Z=$ point , $Y = \bP^m$.\\
Given a point in the moduli space
$$C \stackrel{f}{\rightarrow} Y_l \to Y=\bP^m=\bP^{m-1}\cup \rC^m,$$
pick a point $q=(v_1,v_2, \cdots, v_m) \in \rC^m$ such that $q\notin f(C)$.\\
The one parameter family $\nu :\rC \to PGL(\rC^{m+1})=Aut (\bP^m)$ defined by $\nu (t)= $
$\left( \begin{array}{cccccc}
1 & 0 & 0 & \cdots & 0 & -tv_1\\
0 & 1 & 0 & \cdots & 0 & -tv_2\\
0 & 0 & 1 & \cdots & 0 & -tv_3\\
\vdots & \vdots & \vdots & \vdots & \vdots & \vdots\\
0 & 0 & 0 & \cdots & 1 & -tv_m\\
0 & 0 & 0 & \cdots & 0 & 1\\
\end{array} \right)$
preserves the divisor $\bP^{m-1} \subset \bP^m$.
We use this one parameter family to move $(C, f)$, and note the transformation doesn't change the contact order of $(C, f)$ with the divisor. When $t=1$, $\nu(1)\circ (C,f)$ is supported away from the origin $Z$.
\end{proof}

Proposition~\ref{typeI} and Proposition~\ref{typeII} implies the following:
\begin{theorem}{\label{birationality}}
Suppose each connected component of the submanifold $Z = \coprod_{i} Z_i \subset X$ is of type I or type II. Then we have $\phi_*[\moduli (\tilde{\mathcal{W}_0}, 0, n , \pi^! \beta)]^\vir 
= [\moduli (\mathcal{W}_0, 0, n , \beta)]^\vir$.
\end{theorem}
The following numerical form is a direct consequence of the previous theorem.
\begin{theorem}{\label{numerical}}
Suppose each connected component of the submanifold $Z = \coprod_{i} Z_i \subset X$ is of type I or type II.
Let $V$ be a vector bundle over $X$. Let $\mathbf{c}$ be an invertible multiplicative characteristic class. Then we have an equality of genus-$0$ twisted Gromov-Witten invariants
$$\langle \alpha_1, \cdots, \alpha_n \rangle _{0,n, \beta}^{X, \mathbf{c}, V}=\langle \pi^*\alpha_1, \cdots, \pi^*\alpha_n \rangle _{0,n, \pi^!\beta}^{\tX, \mathbf{c}, \pi^*V}, \text{ where } \alpha_i \in H^*(X) \text{ for all } i.$$
\end{theorem}
\begin{proof}
Since the degeneration used here comes from the deformation to the normal cone from blow-up construction, all insertions involved in the equality, i.e. cohomology classes from $X$ and the vector bundle $V$, can be lifted to the degeneration.
\end{proof}
\subsection{Descendant invariants}

The upshot of this subsection is the following:
\begin{corollary}{\label{descendant}}
Suppose each connected component of the submanifold $Z = \coprod_{i} Z_i \subset X$ is of type I or type II.
If $a_i\leq \text{max}(2, \codim(Z, X)-1)$ for all $i$, then we have
$$\langle \tau_{a_1}\alpha_1, \cdots, \tau_{a_n}\alpha_n \rangle _{0,n, \beta}^X = \langle \tau_{a_1}\pi^*\alpha_1, \cdots, \tau_{a_n}\pi^*\alpha_n \rangle _{0,n, \pi^!\beta}^{\tX}, \text{ where } \alpha_i \in H^*(X).$$
\end{corollary}
If there are too many cotangent line classes $\psi_i$, the previous equality of descendant invariants is not expected to hold. This is because the stabilization of the domain curve via $\vp: \moduli_{0,n}(\tX,\pi^!\beta) \to \moduli_{0,n}(X,\beta)$ causes $\psi_i \neq \vp^*\psi_i$. Indeed, $\psi_i - \vp^*\psi_i$ corresponds to boundary strata in the moduli space.

Given $\tbeta\in H_2(\tX)$, if $\moduli_{0,n}(X,\pi_*\tbeta)$ makes sense, then define
$$ \langle \tau_{a_1}\tau'_{b_1}\gamma_1, \cdots, \tau_{a_n}\tau'_{b_n}\gamma_n \rangle _{0,n, \tbeta}^{\tX} := \int_{[{\moduli_{0,n}(\tX,\tbeta)}]^{\vir}}
\psi_1^{a_1}\cap\vp^*\psi_1^{b_1} \cap\cdots \cap \psi_n^{a_n}\cap\vp^*\psi_n^{b_n}\cap \ev^*(\otimes_{i=1}^n \gamma_i),$$
where $\vp: \moduli_{0,n}(\tX,\tbeta) \to \moduli_{0,n}(X,\pi_*\tbeta)$ and $\gamma_i \in H^*(\tX)$.\\

Theorem~\ref{birationality} implies 
$$\langle \tau'_{b_1}\pi^*\alpha_1, \cdots, \tau'_{b_n}\pi^*\alpha_n \rangle _{0,n, \pi^!\beta}^{\tX}= \langle \tau_{b_1}\alpha_1, \cdots, \tau_{b_n}\alpha_n \rangle _{0,n, \beta}^X.$$
However, $\psi_i \neq \vp^*\psi_i$. In order to prove Corollary~\ref{descendant}, we will show that the correction term vanishes if there are not too many contagent line classes. 

May assume $Z$ is connected. One can blow up successively to deduce results for disconnected submanifold $Z$.
Given an arbitrary map $\pi: Y \to X$, suppose $\pi_*(\beta)=0 \in H_2(X)$, where $\beta \in H_2(Y)$.
Therefore we have a diagram:
$$
\xymatrix{
\moduli_{g,n}(Y, \beta) \ar[d]_p \ar[r]^{\ev} & Y^n \ar[d]^{\pi^n}\\
X \ar[r]^{\Delta} & X^n
}
$$
Suppose $$\Theta \in H^*(\moduli_{g,n}(Y, \beta)), \quad \alpha_i \in H^*(X) \text{ and}\quad \gamma_i \in H^*(Y).$$
For convenience, denote $\moduli_{g,n}(Y, \beta)$ by $\bM$.
\begin{lemma}
We have 
$$\int_{[\bM]^\vir} \Theta\cap \ev^*\big(\otimes_i ( \gamma_i \cap \pi^*\alpha_i ) \big)= 
\int_{p_*\big( [\bM]^\vir \cap \Theta\cap \ev^*(\otimes_i  \gamma_i) \big)}  \cap_i \alpha_i
$$
\end{lemma}
\begin{proof}
This follows from projection formula.
\end{proof}
We will set $Y$ as $\tX$, and $\beta$ as $de$, where $e$ is the $\bP^1$ line class in the exceptional divisor. The previous lemma says we can freely reorganize factors from $H^*(X)$.
\begin{lemma}
Suppose $d>0$. Then $\langle \tau_k \pi^*\alpha , \gamma \rangle_{0,2, de}^{\tX}=0$ when $k \leq \text{max}(1, \codim(Z, X)-2)$.
\end{lemma}
\begin{proof}
By the previous lemma,
$$\langle \tau_k \pi^*\alpha , \gamma \rangle_{0,2, de}^{\tX}= \langle \tau_k  , \pi^*\alpha\cap\gamma \rangle_{0,2, de}^{\tX}.$$
The case $k=0$ is trivial.
When $k=1$, $\langle \tau_1  , \pi^*\alpha\cap\gamma \rangle_{0,2, de}^{\tX}= (2g-2+2)\langle  \pi^*\alpha\cap\gamma \rangle_{0,1, de}^{\tX}=0$.\\

If the invariant doesn't vanish, then we have $\deg (\pi^*\alpha\cap\gamma) \leq \dim X -1$. Otherwise, the pull back of $\pi^*\alpha\cap\gamma$ to the exceptional divisor $D$ is zero. Since $\moduli_{0,2}(\tX,de) \cong \moduli_{0,2}(D,de)$, the invariant vanishes.
On the other hand, $k+ \deg (\pi^*\alpha\cap\gamma) = \vdim \moduli_{0,2}(\tX,de)= \dim X -3 +2 + d (\codim(Z, X) -1)$. Therefore,
$$k \geq d(\codim(Z, X) -1)\geq \codim(Z, X) -1.$$
\end{proof}

Given $\vp: \moduli_{0,n}(\tX,\tbeta) \to \moduli_{0,n}(X,\pi_*\tbeta)$,
the boundary strata associated to $\psi_1 - \vp^*\psi_1$ are given by the clutching morphism from:
$$\moduli_{0,1+\{q\}}(\tX, de) \text{ and } \moduli_{0,\{q'\} + (n-1)}(\tX, \tbeta -de), \text{ where } d \text{ runs through all positive integers.}$$
And then glue two points $q$ and $q'$,
$$\moduli_{0,1+\{q\}}(\tX, de) \times_{\tX} \moduli_{0,\{q'\} + (n-1)}(\tX, \tbeta -de) \to \moduli_{0,n}(\tX, \tbeta).$$
Now we pull back line bundles $\mathbb{L}_1$ and $\vp^*\mathbb{L}_1$ on $\moduli_{0,n}(\tX, \tbeta)$ to $\moduli_{0,1+\{q\}}(\tX, de) \times_{\tX} \moduli_{0,\{q'\} + (n-1)}(\tX, \tbeta -de)$.
$$
\left\{
\begin{array}{ll}
\big(\mathbb{L}_1 \text{ on } \moduli_{0,n}(\tX, \tbeta)\big)|_{\text{strata}} & = \mathbb{L}_1 \text{ on } \moduli_{0,1+\{q\}}(\tX, de).\\
\big(\vp^*\mathbb{L}_1 \text{ on } \moduli_{0,n}(\tX, \tbeta)\big)|_{\text{strata}} & = \theta^*\mathbb{L}_{q'} \text{ on } \moduli_{0,\{q'\} + (n-1)}(\tX, \tbeta -de)\\ &
, \text{ where } \theta : \moduli_{0,\{q'\} + (n-1)}(\tX, \tbeta -de) \to \moduli_{0,\{q'\} + (n-1)}(X, \pi_*\tbeta).
\end{array}
\right.
$$
Suppose $[\Delta]^{\vee}= \sum_i T_i \otimes T_i^\vee$ is the Kunneth decomposition of the Poincare dual of the class $[\Delta]$, where $\Delta: \tX \to \tX \times \tX$ is the diagonal map.
\begin{lemma}
Suppose $\pi: \tX = Bl_Z X \to X$ is an arbitrary blow-up.
If $1\leq a_1\leq \text{max}(2, \codim(Z, X)-1)$, then we have
$$\langle \tau_{a_1}\tau'_{b_1}\pi^*\alpha_1, \tau_{a_2}\tau'_{b_2}\gamma_2, \cdots, \tau_{a_n}\tau'_{b_n}\gamma_n \rangle _{0,n, \tbeta}^{\tX}=
\langle \tau_{a_1+b_1}\pi^*\alpha_1, \tau_{a_2}\tau'_{b_2}\gamma_2, \cdots, \tau_{a_n}\tau'_{b_n}\gamma_n \rangle _{0,n, \tbeta}^{\tX}, 
$$
where $\alpha_1 \in H^*(X)$ and $\gamma_i \in H^*(\tX)$.
\end{lemma}
\begin{proof}
Use the induction on $a_1$.
The analysis of $\psi_1-\vp^*\psi_1$ shows:
\begin{align*}
\langle \tau_{a_1}\tau'_{b_1}\pi^*\alpha_1, \tau_{a_2}\tau'_{b_2}\gamma_2, \cdots, \tau_{a_n}\tau'_{b_n}\gamma_n \rangle _{0,n, \tbeta}^{\tX} =
\langle \tau_{a_1-1}\tau'_{b_1+1}\pi^*\alpha_1, \tau_{a_2}\tau'_{b_2}\gamma_2, \cdots, \tau_{a_n}\tau'_{b_n}\gamma_n \rangle _{0,n, \tbeta}^{\tX}\\
+ \sum_{d>0} \sum_{T_i} \langle \tau_{a_1-1}\pi^*\alpha_1, T_i \rangle _{0,2, de}^{\tX} \bullet
\langle \tau'_{b_1} T_i^{\vee}, \tau_{a_2}\tau'_{b_2}\gamma_2, \cdots, \tau_{a_n}\tau'_{b_n}\gamma_n \rangle _{0,n, \tbeta -de}^{\tX}.
\end{align*}
Since $a_1-1 \leq \text{max}(1, \codim (Z, X)-2)$,  by the previous lemma $\langle \tau_{a_1-1}\pi^*\alpha_1, T_i \rangle _{0,2, de}^{\tX} =0$.
\end{proof}
\begin{proof}[Proof of Corollary~\ref{descendant}]
In the previous lemma, set $\gamma_i = \pi^*\alpha_i$, $\tbeta=\pi^!\beta$ and $b_1=b_2=\cdots=b_n=0$.
Then apply the lemma to $a_1, a_2, \cdots, a_n$, this shows
$$
\langle \tau_{a_1}\pi^*\alpha_1, \cdots, \tau_{a_n}\pi^*\alpha_n \rangle _{0,n, \pi^!\beta}^{\tX}=
\langle \tau'_{a_1}\pi^*\alpha_1, \cdots, \tau'_{a_n}\pi^*\alpha_n \rangle _{0,n, \pi^!\beta}^{\tX}.
$$
Now it follows from Theorem~\ref{birationality}.
\end{proof}
\subsection{Examples and remarks}

\begin{example}{\label{source}}
Given any projective manifold $X$, here we provide several ways to find a submanifold $Z\subset X$, so that $N_{Z/X}$ is generated by global sections.
\begin{enumerate}
\item Embed $X$ in a homogeneous space $\mathcal{P}$, and pick an arbitrary submanifold $Y \subset \mathcal{P}$. By Bertini's Theorem, one can perturb $Y$ so that $Y$ is transversal to $X$. Then take $Z= X \cap Y$.
\item Take any vector bundle $V$ over $X$ and an ample line bundle $L$. When $n >> 0$, $V\otimes L^{\otimes n}$ is generated by global sections. Take a regular section $s\in H^0(X, V\otimes L^{\otimes n})$, and let $Z =s^{-1}(0)$.
\item Suppose $L_1, L_2, \cdots, L_m$ are line bundles over $X$, and are generated by global sections. Take a regular section $s\in H^0(X, \oplus_{i=1}^m L_i)$, and let $Z =s^{-1}(0)$. Then $Z$ is a complete intersection of $X$.
\end{enumerate}
\end{example}

\bigskip

\begin{example}
Suppose $L_1, L_2, \cdots, L_m$ are line bundles over $Z$, and each $L_i$ is generated by global sections.
Let $X=\bP_Z (\oplus_{i=1}^m L_i \oplus \sO_Z)$ be the projective completion, and $\tX$ be the blow-up along $Z$. Let $(\rC^*)^m$ act on $X$ and $\tX$ fiberwisely. In principle, one can use virtual localization to express all GW-invariants of $\tX$ and $X$ in terms of those of $Z$, and use the calculation to prove Theorem~\ref{numerical} in this case. However, this is a formidable combinatorial task. When $Z$ is a point and $\beta\in H_2(X)$ has degree 2, Theorem~\ref{numerical} corresponds to Lagrangian interpolation in the localization computation after cancelling numerous terms.
\end{example}

\bigskip

\begin{remark}
Suppose $N_{Z/X}$ is generated by global sections, and $\pi: \tX \to X$ is the blow-up.
Let $V$ be a convex line bundle over $X$, and $s\in H^0(X,V)$ is a section so that $Y:= s^{-1}(0)$ is a submanifold of $X$. It is well-known that genus-$0$ GW-invariants of $Y$ can be expressed by twisted invariants of $X$. If $Y$ and $Z$ is transversal in $X$, then $\pi^*(s) \in H^0(\tX,\pi^*V)$ is a regular section. And we have a Cartesian diagram:
$$
\xymatrix{
\tY=Bl_{Y\cap Z} Y = \pi^*(s)^{-1}(0) \ar[d] \ar[r] & \tX \ar[d]\\
Y= s^{-1}(0) \ar[r] & X.
}
$$
Since $\pi^*V$ is also a convex line bundle of $\tX$, by Theorem~\ref{numerical}, we have
$$
\{GW(\tY) \text{ with insertions from } Y \}=  \{\text{ twisted- }GW(\tX) \text{ with insertions from } X \}
=  \{\text{ twisted- }GW(X)\} = \{GW(Y) \}.
$$
On the other hand, $N_{(Y\cap Z) / Y}$ is the pull back of $N_{Z/X}$, and is generated by global sections as well. This also implies 
$\{GW(\tY) \text{ with insertions from } Y \}=\{GW(Y) \}$.

For arbitrary blow-ups, the correction terms between GW-invariants of $\tX$ and $X$ are required. If the universal blow-up formula exists, the correction terms should also have this functoriality.
\end{remark}

\bigskip

\begin{remark}{\label{difficulty1}}
We speculate that Theorem~\ref{numerical} holds as long as $N_{Z/X}$ is convex without any additional assumption. Here we briefly discuss the technical difficulty encountered in our approach. For simplicity, assume $X = \bP_Z(N \oplus \sO)$. Given any point $(C, f)\in \moduli_{0,n}(X, \beta)$, we have $C \stackrel{f}{\rightarrow} X \stackrel{p}{\rightarrow} Z$. Because $N$ is convex, $f^*p^*N$ is generated by global sections. Therefore one can find a section of $f^*p^*N$ to move the curve so that the new curve is supported away from $Z$ in $X$.

On the other hand, suppose $\C_i$ is an irreducible component of the virtual normal cone, and $(C,f) \in \text{supp} (\C_i)$. To check the technical assumption of Proposition~\ref{special}, one has to make sure that the new curve still stays in $\text{supp} (\C_i)$. A priori, $\C_i$ can be supported in a very small part of $\moduli_{0,n}(X,\beta)$ (even if one assumes the reduced structure of $\moduli_{0,n}(X,\beta)$ is smooth). More precisely, $\dim \supp(\C_i) \geq \vdim \moduli_{0,n}(X,\beta)$, and the equality can be achieved. Local analysis of Kuranish map is required if one attempts to prove in this way.
\end{remark}

\bigskip

\begin{example}{\label{counter1}}
In this example we will see that even if $Z \subset X$ has freedom to move to avoid any finite collection of holomorphic curves in $X$, the induced GW-invariants of $\tX$ can be different from the corresponding GW-invariants of $X$.

Consider two vector bundles on $Z=\bP^r$ : $ V_1= \oplus_{i=1}^s \sO_Z$  and $ V_2= \oplus_{i=1}^t \sO_Z(-k)$, where $s,t \geq 2$ and $k>0$. Let $X$ be the projective completion $\bP_Z (V_1 \oplus V_2 \oplus \sO_{\infty})$, and $Z \subset X$ with normal bundle $V_1\oplus V_2$. Since $s\geq 2$, $Z$ can be moved to avoid any finite collections of holomorphic curves. Define
$$
\begin{array}{l}
\pi : \tX \to X \text{ is the blow-up along } Z, \\
Y := \bP_Z (V_1 \oplus \sO_{\infty}) \subset X, \\
\pi_Y : \tY \to Y \text{ is the blow-up along } Z\subset Y.
\end{array}
$$
We have a diagram (not Cartesian):
$
\xymatrix{
\tY = Z \times Bl_{\{0\}} \bP^s \ar[d] \ar[r] & \tX \ar[d] \\
Y= Z \times \bP^s \ar[r] & X.
}
$
\qquad
$
\xymatrix{
\moduli_{0,n+1}(Z, d\ell) \ar[r] \ar[d]_{\pi_{n+1}} & Z\\
\moduli_{0,n}(Z, d\ell).
}
$\\
Let $[\ell]$ be the line class in $Z$. Define an obstruction bundle on $\moduli_{0,n}(Z, d\ell)$ by $\calU_d := R^1\pi_{{n+1}_*}\ev^* \sO_Z (-k)$.\\
Regard $\ell$ as a curve class in $X$ via $Z \subset X$. Let $\Phi : \moduli_{0,n}(\tX, \pi^![d\ell]) \to \moduli_{0,n}(X, [d\ell])$ with $d>0$.

\begin{lemma}{\label{example}}
\mbox{}
\begin{enumerate}
\item $\moduli_{0,n}(X, d\ell) \cong \moduli_{0,n}(Z, d\ell) \times \bP^s$,
\item $\moduli_{0,n}(\tX, \pi^!d\ell) \cong \moduli_{0,n}(Z, d\ell) \times Bl_{\{0\}}\bP^s$,
\item $ Ob\big( \moduli_{0,n}(X, d\ell) \big) \cong \calU_d \boxtimes \big( \oplus_t \sO_{\bP^s}(1) \big)$,
\item $ Ob\big( \moduli_{0,n}(\tX, \pi^!d\ell) \big) \cong \calU_d \boxtimes \big[  \oplus_t \big( \vp^*\sO_{\bP^s}(1) \otimes [-D]\big) \big]$, where $\vp : Bl_{\{0\}}\bP^s \to \bP^s$ and $D$ is the exceptional divisor of $\vp$.
\end{enumerate}
\end{lemma}
Note $\Phi$ is a birational map in this case, but the natural map between obstruction bundles is not surjective.
Assume the lemma, then the difference of (push-down) virtual classes $\Phi_*[\moduli_{0,n}(\tX, \pi^!d\ell)]^\vir - [\moduli_{0,n}(X, d\ell)]^\vir$ in general doesn't vanish, and has non-zero contribution to GW-invariants. For example, take $s=t=k=2$ and $d=1$, then $\calU_d$ is a trivial line bundle. Let $H$ be the hyperplane class of $\bP^2$. Then 
\begin{align*}
& \vp_* [(H-D)^2] - H^2= -[pt] \in A_0(\bP^2)\\
\Longrightarrow & \Phi_*[\moduli_{0,n}(\tX, \pi^!\ell)]^\vir - [\moduli_{0,n}(X, \ell)]^\vir = -[\moduli_{0,n}(Z, \ell)]^\vir \times \{pt\} \in A_*(\moduli_{0,n}(Z, \ell)\times \bP^2),\\
& \text{ which apparently has non-zero contribution to GW-invariants. }
\end{align*}
For general $s, t, k, d$, the difference of (push-down) virtual classes is given by $[\text{ twisted -}\moduli_{0,n}(Z, d\ell)]^\vir \times \{pt\}$, where the virtual class is twisted by the vector bundle $V_2 \to Z$, and the characteristic class is a combination of various chern classes.

Now we sketch the proof of Lemma~\ref{example}. First note the normal bundle $N_{Y/X} \cong \sO_Z(-k)\boxtimes \big( \oplus_t \sO_{\bP^s}(1) \big)$. This vector bundle is $[d\ell]$-concave (but is not concave for any curve class), therefore $\moduli_{0,n}(Y,d\ell)$ is a path-connected component of $\moduli_{0,n}(X,d\ell)$. Let $[\ell_1]= [\ell]$ and $[\ell_2]$ be the line classes of $Z$ and $\bP^s$.
The equality $\moduli_{0,n}(X,d\ell) = \moduli_{0,n}(Y,d\ell)$ follows from the following lemma.
\begin{lemma}{\label{factor1}}
For any $f; \bP^1 \to X$, if $f(\bP^1)\nsubseteq Y$, then 
$$f_*[\bP^1] = a [\ell_1]+b[\ell_2] \in A_1(X) \cong A_1(Y)=A_1(Z)\oplus A_1(\bP^s) \text{ with } a\geq 0, b>0.$$
\end{lemma}
The obstruction sheaf on $\moduli_{0,n}(Y,d\ell)$ is deduced directly from the normal bundle $N_{Y/X}$.

\begin{lemma}
Given three manifolds $Z \subset Y \subset X$, we have a diagram (not Cartesian):
$$
\xymatrix{
\tY = Bl_Z Y \ar[d] \ar[r] & \tX= Bl_Z X \ar[d] \\
Y \ar[r] & X.
}
$$
Then $N_{\tY/\tX} \cong \pi^*(N_{Y/X})\otimes [-D]$, where $D$ is the exceptional divisor of $\tY \to Y$.
\end{lemma}

In our case, the lemma says $N_{\tY/\tX}\cong \sO_Z(-k)\boxtimes \big[  \oplus_t \big( \vp^*\sO_{\bP^s}(1) \otimes [-D]\big) \big]$, which is also a $[d\ell]$-concave bundle. An analogue of Lemma~\ref{factor1} shows  $\moduli_{0,n}(\tY,d\ell)$ is the only component of $\moduli_{0,n}(\tX,\pi^!d\ell)$. The analysis of obstruction sheaf is straightforward.
\end{example}
\section{Vanishing Theorems for blow-ups}
Suppose we have a map $f: X \to Y$ between two compact complex manifolds. It is obvious that 
$$ \int_X{\alpha \wedge f^*\beta}=0, \alpha \in H^*(X), \beta \in H^*(Y), \text{ when } \deg_{\mathbb{R}}\beta > 2\dim_{\mathbb{C}} Y.$$
However, the virtual version in general is not true (even if $X$ and $Y$ are smooth):
$$ \int_{[X]^\vir}{\alpha \wedge f^*\beta}\stackrel{?}{=}0, \alpha \in H^*(X), \beta \in H^*(Y), \text{ when }\deg_{\mathbb{R}}\beta > 2\vdim Y.$$
To rectify this situation, one has to impose the assumption that
\textbf{
$f: X \to Y$ has compatible perfect obstruction theories.}
With such assumption, the vanishing result holds in the virtual version. This simple phenomenon is the starting point of vanishing theorems for blow-up in this paper.\\

In our convention, the empty set $\emptyset$ has dimension $= -\infty$, and $\codim(\emptyset, S)= +\infty$ if $S$ is not empty. 

\begin{lemma}{\label{vanish}}
Suppose $\bM$ and $\bN$ are two proper DM-stacks. Assume $\vp : \bM \to \bN$ has compatible perfect obstruction theories. Let $\alpha \in A^*(\bM), \beta \in A^*(\bN)$. Denote the virtual normal cone of $\bN$ by $\C_{\bN}$.
If there exists an open substack $U \subset \bN$ such that 
\begin{enumerate}
\item $\vp(\bM)\cap U= \emptyset$
\item $\dim (\C_{\bN}-\C_{\bN}|_U) \leq \dim \C_{\bN} - k$
\end{enumerate}
Then $\int_{[\bM]^\vir} \alpha \cap \vp^*\beta=0$ when deg$\beta > \vdim \bN - k$.
\end{lemma}
\begin{proof}
We will adapt notation used in Section 3. Let $\Ldot$, $\Edot$ and $\Fdot$ be the compatible perfect obstruction theories on $\bM/\bN$, $\bM$ and $\bN$ respectively. Recall the diagram
$$
\xymatrix{
\C_{\bN} & \ar[l]_(.7){\psi}\C_{\bM} \cap \vectvF \ar[r] \ar[d]& \C_{\bM} \ar[d]\\
& \vectvF \ar[r]^{0_{\sL_2}} & \vectLvF 
}
$$
where $\C_{\bM}$ and $\C_{\bN}$ are virtual normal cones used to construct virtual classes. Note $\psi$ is a proper map.
By abuse of notation, $\alpha$ (and $\beta$) will be also viewed as an element in $ A^*(\C_{\bM})$ (and $A^*(\C_{\bN})$).
\begin{align*}
\vp_* \big( [\bM]^{\vir}\cap \alpha\cap\vp^*\beta \big) &= 
\vp_* \circ 0^!_{(\vp^*\F_2)} \circ 0_{\sL_2}^![\C_{\bM}\cap\alpha\cap\psi^*\beta] = 
0^!_{\F_2} \circ \psi_* \circ 0_{\sL_2}^![\C_{\bM}\cap\alpha\cap\psi^*\beta]\\&= 
0^!_{\F_2} \big( \beta\cap(\psi_* \circ 0_{\sL_2}^![\C_{\bM}\cap\alpha]) \big) , \text{ where }
\psi_* \circ 0_{\sL_2}^![\C_{\bM}\cap\alpha] \in A_*\big(\psi(\C_{\bM} \cap \vectvF )\big).
\end{align*}
Since $\vp(\bM)\cap U= \emptyset$, we have $\psi\big(\C_{\bM} \cap \vectvF \big) \subset \C_{\bN}-\C_{\bN}|_U$.
Recall that $\dim \C_{\bN}=\vdim (\bN)+\rk(\F_2)$. Therefore $\dim \psi\big(\C_{\bM} \cap \vectvF \big) \leq \vdim (\bN)+\rk(\F_2) -k$.
Because $\deg \beta + \rk(\F_2) > \vdim \bN + \rk(\F_2) -k$, we know $0^!_{\F_2} \circ (\beta \cap *) : A_{\bullet}\big(\psi (\C_{\bM} \cap \vectvF )\big) \to A_{\bullet -\text{ deg}\beta - \rk(\F_2)}(\bN)$ must be a zero map.
\end{proof}
There is a topological statement of Lemma~\ref{vanish}. All $A_*(\bullet)$ in the proof must be replaced by Borel-Moore homology $H_*^{BM}(\bullet)$. The proof is the same and is omitted.
\begin{lemma}{\label{topovanish}}
Suppose $\bM$ and $\bN$ are two proper DM-stacks. Assume $\vp : \bM \to \bN$ has compatible perfect obstruction theories. Let $\alpha \in H^*(\bM), \beta \in H^*(\bN)$. 
If there exists an open substack $U \subset \bN$ such that 
\begin{enumerate}
\item $\vp(\bM)\cap U= \emptyset$
\item $\dim (\C_{\bN}-\C_{\bN}|_U) \leq \dim \C_{\bN} - k$
\end{enumerate}
Then $\int_{[\bM]^\vir} \alpha \cap \vp^*\beta=0$ when $\deg_{\mathbb{R}}\beta > 2\vdim \bN - 2k$.
\end{lemma}
\begin{remark}
\mbox{}
\begin{enumerate}
\item The second assumption $\dim (\C_{\bN}-\C_{\bN}|_U) \leq \dim \C_{\bN} - k$ only depends on $U$ and the singularities of $\bN$, but is independent of the perfect obstruction theory $\Fup$.
\item Taking $U$ as an empty set and $k=0$, this is the vanishing result mentioned at the beginning of this section.
\end{enumerate}
\end{remark}
\begin{corollary}{\label{vanish2}}
Suppose $\bM$ and $\bN$ are two proper DM-stacks. Assume $\vp : \bM \to \bN$ has compatible perfect obstruction theories.
Suppose $A \subset B$ is a pair of compact complex manifolds, with a fiber diagram:
$$
\xymatrix{
\bM' \ar[d]^{j'} \ar[r]^{i'} & \bM \ar[d]_j \ar[r]^\vp & \bN \\
A \ar[r]^i & B.
}
$$
If there exists an open substack $U \subset \bN$ such that 
\begin{enumerate}
\item $\vp\circ i'(\bM')\cap U= \emptyset$
\item $\dim (\C_{\bN}-\C_{\bN}|_U) \leq \dim \C_{\bN} - k$
\end{enumerate}
Then $\int_{[\bM]^\vir} j^*\big(PD_B\circ i_*(w)\big)   \cap  \alpha \cap \vp^*\beta=0$ when $\deg_{\mathbb{R}}\beta > 2\vdim \bN - 2k$.\\
Here $\alpha \in H^*(\bM), \beta \in H^*(\bN)$, $w\in H_*(A)$ and $PD_B$ is the Poincare dual in $B$.
\end{corollary}

\begin{proof}
Form a fiber diagram:
$$
\xymatrix{
\C_{\bM}' \ar[r]^{i''} \ar[d]^{p'} & \C_{\bM} \ar[d]^{p} \\
\bM' \ar[d]^{j'} \ar[r]^{i'} & \bM \ar[d]_j \\
A \ar[r]^i & B.
}
$$
$$
[\C_{\bM}] \cap p^*j^*\big(PD_B\cdot i_*(w)\big) =  i''_*\big(p'^* j'^* PD_A(w) \cap i^![\C_{\bM}] \big) \in H_*^{BM}(\C_{\bM}).
$$
Here $i^!$ means cap with $j^*(u_{A,B})$, where $u_{A,B} \in H^*(B, B-A)$ is the canonical orientation class of $A\subset B$. Note $i^![\C_{\bM}] \in H_*^{BM}(\C_{\bM}')$. Let $\gamma:= p'^* j'^* PD_A(w) \cap i^![\C_{\bM}]$ and $\C_{\bM}'\cap \vect(i'^*\vp^*\calF_2) \stackrel{\psi'}{\rightarrow} \C_{\bN}$. Note $\psi'$ is a proper map. We have 
$$\int_{[\bM]^\vir} j^*\big(PD_B\circ i_*(w)\big)   \cap  \alpha \cap \vp^*\beta = 0^!_{\F_2} \big( \beta\cap(\psi'_* \circ 0_{\sL_2}^![\gamma\cap\alpha]) \big).$$
Now argue as the proof of Lemma~\ref{vanish} and note the image of $\psi'$ lies in $\C_{\bN}-\C_{\bN}|_U$.
\end{proof}

\subsection{Relative case}
In this subsection, we always assume $Z$ is connected.
Suppose $N_{Z/X}$ is a convex bundle.
We will attempt to apply the vanishing lemma to 
$$\vp : \bM =\moduli_{0,n}(\tX,\tbeta) \to \moduli_{0,n}(X, \pi_*\tbeta) \to \moduli_{0,m}(X, \pi_*\tbeta)=\bN, \tbeta \neq \pi^!\pi_*\tbeta,$$ 
where the second arrow forgets the last $n-m$ marked points. Note
$$ \tbeta \neq \pi^!\pi_*\tbeta \Longleftrightarrow \tbeta = \pi^!\beta + de , d \neq 0 ,\text{ where } e \text{ is the line class in the exceptional divisor} .$$
The open substack $U \subset \bN$ will be a collection of stable maps supported away from the submanifold $Z \subset X$. To show the composition map $\vp$ has compatible perfect obstruction theories, note the first map has compatible perfect obstruction theories (Proposition~\ref{compatible}), and so does the forgetful map.

Unfortunately, it is difficult to directly check the second assumption in Lemma~\ref{topovanish} if $k>0$. Degeneration formula will be used to simplify the situation.

First we consider the simplest case: $Z=$ the origin $\subset X = \bP^r$, with the divisor $D = \bP^{r-1} \subset X$.
Let $\bN = \moduli (\bP^r, \bP^{r-1} , \Gamma)$, where $\Gamma$ is an adimissible graph. In this case, $\bN$ is a smooth DM-stack. 
\begin{lemma}
We have $\codim (\bN-U,\bN) \geq r-1$.
\end{lemma}
\begin{proof}
Define $\nu :\rC^r \to PGL(\rC^{r+1})=Aut (\bP^m)$ by 
$\nu (v_1,v_2,\cdots,v_r)= $
$\left( \begin{array}{cccccc}
1 & 0 & 0 & \cdots & 0 & -v_1\\
0 & 1 & 0 & \cdots & 0 & -v_2\\
0 & 0 & 1 & \cdots & 0 & -v_3\\
\vdots & \vdots & \vdots & \vdots & \vdots & \vdots\\
0 & 0 & 0 & \cdots & 1 & -v_r\\
0 & 0 & 0 & \cdots & 0 & 1\\
\end{array} \right)$
This matrix preserves the divisor $\bP^{r-1}$ and doesn't not change the contact order of the curve to $\bP^{r-1}$, and therefore induces an action on $\bN = \moduli (\bP^r, \bP^{r-1} , \Gamma)$. 

Equip $\bN -U$ with reduced structure, denote it by $B$. Suppose $\codim (B,\bN) < r-1$, then there exists a point $[C \stackrel{f}{\rightarrow} Y_{l} \to Y] \in B$, so that $B$ is smooth at the point $[C \stackrel{f}{\rightarrow} Y_{l} \to Y]$ and $\codim (B,\bN) < r-1$ near the point $[C \stackrel{f}{\rightarrow} Y_{l} \to Y]$.
Define a morphism $\sigma:\rC^r \to \bN$ by the action of $\rC^r$ on $[C \stackrel{f}{\rightarrow} Y_{l} \to Y]$.
$$\sigma(v_1,\cdots,v_r) \in B \Longleftrightarrow (v_1,\cdots,v_r) \in \text{ the image of the curve } C.$$
Therefore $\dim \sigma^{-1}(B) \leq 1$. Take the linearized map of $\sigma$:
$$ T\sigma|_Z : T\rC^r|_Z \to T\bN|_{[C \stackrel{f}{\rightarrow} Y_{l} \to Y]}.$$
$$ \dim T\sigma|_Z^{-1}\big(TB|_{[C \stackrel{f}{\rightarrow} Y_{l} \to Y]}\big) \leq 1.$$
Therefore $\codim (B,\bN) \geq \codim\{ T\sigma|_Z^{-1}\big(TB|_{[C \stackrel{f}{\rightarrow} Y_{l} \to Y]}\big) , T\rC^r|_Z\} \geq r-1$, which is a contradiction.
\end{proof}
Now we can control the codimension for type II cases.
\begin{lemma}{\label{codimtypeII}}
Suppose $Z$ is of type II. Let $Y= \bP(N_{Z/X}\oplus \sO)$, $D = \bP(N_{Z/X})$ and $\bN = \moduli_{0,n}(Y, D, \Gamma)$.
Then $$\codim(\C_{\bN}-\C_{\bN}|_U, \C_{\bN}) \geq \rk(N_{Z/X}) -1.$$
\end{lemma}
\begin{proof}
$\bN$ is a smooth DM-stack. It suffices to prove $\codim(\bN-U, \bN) \geq \rk(N_{Z/X}) -1$.\\
The fibration $\bN = \moduli_{0,n}(Y, D, \Gamma) \to Z$ is locally trivial, and therefore reduces the problem to the fiber.
Now it follows from the previous lemma.
\end{proof}
\begin{lemma}{\label{codimtypeI}}
Suppose $N_{Z/X}$ is convex and there is a subbundle $\mathcal{F} \subset N_{Z/X}$ generated by global sections.
Let $Y= \bP(N_{Z/X}\oplus \sO)$, $D = \bP(N_{Z/X})$ and $\bN = \moduli_{0,n}(Y, D, \Gamma)$.
Then there exists a section $q\in H^0(Z,N_{Z/X})$, so that 
$$\codim(\C_{\bN}-\C_{\bN}|_{U_q}, \C_{\bN}) \geq \rk(\mathcal{F})-1,$$ where $U_q \subset \bN$ is a collection of relative stable maps supported away from $q(Z)$ in $Y$.
\end{lemma}
In particular, if $N_{Z/X}$ is generated by global sections, then we have a good bound of the codimension $\geq \rk(N_{Z/X}) -1$. Now the goal is to prove Lemma~\ref{codimtypeI}. $(\bullet)^\red$ means the reduced structure.
\begin{lemma}
Suppose $f: \calA \to \mathcal{B}$ is a morphism of separated DM-stacks of finite type over $\rC$. Then $\calA$ can be splitted as finite disjoint union $\calA = \coprod_\finite \calA_i$, so that

\begin{enumerate}
\item For each $i$, $\calA_i$ is irreducible and locally closed in $\calA$ and then equipped with reduced structure.
\item Set $f_i : \calA_i \to \overline{f(\calA_i)}$. Then $f_i^{-1}\big(f_i(a)\big)$ has $\dim \calA_i - \dim \overline{f(\calA_i)}$ for all $a\in \calA_i$.
\end{enumerate}
\end{lemma}
\begin{proof}
Use the induction on the number of irreducible components of top dimension in the domain. Suppose $D$ is an irreducible component of top dimension in $\calA$. The induced map $f: D^\red \to \overline{f(D)}^\red$ is a dominant morphism of integral DM-stack of finite type over $\rC$. There exists an open substack $\calU \subset D^\red$, such that for any $y\in f(D^\red)$, $\dim \calU_y = \dim D - \dim \overline{f(D)}$. It remains to consider $f : \calA-D \to \mathcal{B}$ and $f: D^\red-\calU \to \mathcal{B}$.
\end{proof}

Given an adimissible graph $\Gamma$ for $(Y, D)$ and assume $\moduli(Y,D,\Gamma)$ exists. Define 
$$\moduli(Y, \Gamma):= \prod_{v\in V(\Gamma)} \moduli_{g(v), \#\text{legs}(v)+\#\text{roots}(v)} (Y, b(v)).$$
Because $\Gamma$ is relatively connected, $\moduli(Y, \Gamma)$ makes sense and is the moduli space of (disconnected)-stable maps in $Y$. Note here we have used the condition: if $|V(\Gamma)|>1$, then each vertex $v\in V(\Gamma)$ has at least one root and $b(v)\neq 0$.

There is a natural map
$$\moduli(Y,D,\Gamma) \to \moduli(Y, \Gamma).$$
But there is \textbf{no} natural arrow between two obstruction theories.
The universal curve of $\moduli(Y, \Gamma)$ is 
$$\moduli(Y, \Gamma)^\univ = \coprod_{v\in V(\Gamma)} \frac{\moduli(Y, \Gamma)}{\moduli_{g(v), \#\text{legs}(v)+\#\text{roots}(v)} (Y, b(v))} \times \moduli_{g(v), \#\text{legs}(v)+\#\text{roots}(v)+1} (Y, b(v)).$$
Note the coarse moduli space of $\moduli(Y, \Gamma)$ is projective, as shown in~\cite{FP}.

\begin{lemma}
Suppose $\calA$ is a separated DM-stack of finite type over $\rC$ with pure dimension. Assume $\calF \subset N_{Z/X}$ is a subbundle generated by global sections.\\
Given $\vp : \calA \to \moduli(Y, \Gamma)$, then there always exists a section $q\in H^0(Z, N_{Z/X})$,\\ such that $\codim (\calA-\calA|_{\calU_q}, \calA) \geq \rk(\calF)-1$. Here $\calU_q \subset \moduli(Y, \Gamma)$ is defined by (disconnected-) stable maps supported away from $q(Z)$ in $Y$.
\end{lemma}

\begin{proof}
Recall $Y= \bP(N_{Z/X}\oplus \sO)$ and $D = \bP(N_{Z/X})$. Let $r=\rk(\calF)$. May assume $r\geq 2$, otherwise it is trivial.\\
\textbf{\underline {Step 1}}\\
Use the previous lemma to split $\calA = \coprod_\finite \calA_i$. Define $\moduli_i:= \overline{\vp(\calA_i)}^\red$, and $\coarse_i$ is the image of $\moduli_i$ in the coarse moduli space $\coarse(Y, \Gamma)$. Assume $\coarse(Y, \Gamma) \hookrightarrow \bP^N$ and $\dim \moduli_i=\dim \coarse_i= k_i$.\\

If $k_i \geq r-2$, then pick a subplane $\bP^{N-k_i+(r-2)}$ in $\bP^N$ such that $\dim (\bP^{N-k_i+(r-2)}\cap \coarse_i ) = r-2$. Define two new objects by the fiber diagrams:
$$
\xymatrix{
\calW_i^\univ \ar[r] \ar[d] & \moduli(Y, \Gamma)^\univ \ar[d]\\
\calW_i \ar[r] \ar[d] & \moduli(Y, \Gamma) \ar[d]\\
\bP^{N-k_i+(r-2)}\cap \coarse_i \ar[r] & \coarse(Y, \Gamma)
}
$$
If $k_i < r-2$, then define $\calW_i$ as $\moduli_i$.
Note $\dim \calW_i^\univ \leq r-1$ and $f_i : \calW_i^\univ \to Y$ is a proper map. Suppose $\oplus_s \sO_Z \to \calF \to 0$.
$$\rC^s \stackrel{p \text{ projection}}{\longleftarrow} Z\times \rC^s = \vect (\oplus_s \sO_Z) \stackrel{\theta \text{ smooth}}{\longrightarrow} \vect (\calF) \stackrel{closed}{\hookrightarrow} \vect (N_{Z/X}) \stackrel{open}{\hookrightarrow} Y.$$
\begin{align*}
\dim p_*\theta^{-1}\big(f_i(\calW_i^\univ)\cap \vect(\calF)\big) & \leq \dim \theta^{-1}\big(f_i(\calW_i^\univ)\cap \vect(\calF)\big)\\
& = s-\rk{\calF} + \dim \big(f_i(\calW_i^\univ)\cap \vect(\calF)\big)\\
& \leq s-r+r-1=s-1.
\end{align*}
Since there are finitely many $\calW_i^\univ$, there exists $q\in \rC^s$ such that
$$q\notin \bigcup_i p_*\theta^{-1}\big(f_i(\calW_i^\univ)\cap \vect(\calF)\big).$$
Such $q$ induces an section in $H^0(Z, \calF) \subset H^0(Z, N_{Z/X})$.\\
\textbf{\underline {Step 2}}  \qquad Claim $\dim (\moduli_i - \moduli_i|_{\calU_q}) \leq k_i - (r-1)$. \\
Suppose $k_i \geq r-2$.
Let $\calU_q^\comp := \moduli(Y, \Gamma) - \calU_q$ equipped with reduced structure, and $U_q$ be the corresponding coarse moduli.
Argue by contradiction. Suppose not, then $\dim \calU_q^\comp \cap \moduli_i \geq k_i-(r-2)$.
\begin{align*}
& \Longrightarrow \dim U_q^\comp \cap \coarse_i \geq k_i-(r-2).\\
& \Longrightarrow \dim U_q^\comp \cap \coarse_i \cap \bP^{N-k_i+(r-2)} \neq \emptyset \text{ in } \bP^N.\\
& \Longrightarrow \calU_q^\comp \cap \calW_i \neq \emptyset \text{ in } \moduli(Y,\Gamma).
\end{align*}
On the other hand,
$$ q(Z) \cap f_i(\calW_i^\univ)=\emptyset \Longrightarrow \calW_i \subset \calU_q,$$
which is a contradiction.\\
If $k_i< r-2$, then a similar argument shows $\moduli_i \subset \calU_q$. Therefore $\moduli_i - \moduli_i|_{\calU_q}= \emptyset$.\\
\textbf{\underline {Step 3}}\\
$$ \dim (\moduli_i - \moduli_i|_{\calU_q}) \leq \dim \moduli_i - (r-1) \text { and } \vp_i: \calA_i \to \moduli_i \text{ has the fiber dimension } \dim \calA_i-\dim \moduli_i,$$
\begin{align*}
\Longrightarrow \dim (\calA_i - \calA_i|_{\calU_q}) & \leq \dim \moduli_i -(r-1) + \dim \calA_i - \dim \moduli_i\\
& =\dim \calA_i -(r-1) \leq \dim \calA -(r-1).
\end{align*}
Now it follows from $\calA - \calA|_{\calU_q} = \coprod_{\finite} (\calA_i - \calA_i|_{\calU_q})$.
\end{proof}
\begin{proof}[Proof of Lemma~\ref{codimtypeI}]
Consider the composition $\calA= \C_{\moduli(Y,D,\Gamma)} \to \moduli(Y,D,\Gamma) \to \moduli(Y, \Gamma)$.\\
Let $\calU_{\moduli(Y,D,\Gamma),q}$ ( and $\calU_{\moduli(Y,\Gamma),q}$) be the (relative) stable maps supported away from $q(Z)$ in $Y$.\\
Note $\calU_{\moduli(Y,D,\Gamma),q}=$ the preimage of $\calU_{\moduli(Y,\Gamma),q}$ under the natural map. Now it follows from the previous lemma.
\end{proof}

The next two corollaries are the building blocks of vanishing theorems of absolute GW-invariants.
According to Lemma~\ref{codimtypeI} and Lemma~\ref{codimtypeII}, we define the codimension $\delta$ corresponding to $k$ in Lemma~\ref{topovanish} as
$$
\delta= \left\{
\begin{array}{ll}
\rk(\calF) -1 & \text{ if } Z \subset \text{ X is of type I, and } \calF\subset N_{Z/X} \text{ is generated by global sections.}\\
\rk(N_{Z/X})-1 & \text{ if } Z \subset \text{ X is of type II.}
\end{array}
\right.
$$
Recall $\pi : \tY:=\overline{\bP(N_{Z/X} \oplus \sO_Z)} \to Y:=\bP(N_{Z/X} \oplus \sO_Z)$ is the blow up along $Z$, and $D:= \bP(N_{Z/X}) \subset Y$. Let $A\subset [n]$ , $\alpha_a \in H^*(Y)$ for $a\in A$, $\gamma_i\in H^*(\tY)$ for $i\in [n]$, and $t_{\ast} \in H^*(D)$. See the paragraph before Theorem~\ref{thm3} about notation of GW-invariants. 

\begin{corollary}{\label{relvanish1}}
Suppose $Z \subset X$ is of type I or II. Let $\tGamma$ be an admissible weighted graph for $(\tY, D)$, and $\pi_{A*}$ be a composition of push-forward and the map forgetting the $[n]-A$ legs. Assume $\moduli(Y, D, \pi_{A*} \tGamma)$ makes sense.

If $\tGamma$ satisfies both conditions
$$
\left\{
\begin{array}{l}
\text{ genus-zero weight } g: V(\tGamma) \stackrel{\equiv 0}{\rightarrow} \mathbb{Z}_{\geq 0},\\
\text{ homology weight }b(v) \neq \pi^!\pi_*b(v) \text{ for at least one vertex } v \in V(\tGamma).
\end{array}
\right.
$$

Then we have 
$$\langle \longvec{\pi^* \alpha _A} \cdot \longvec{\tau_\bullet \gamma_{[n]}}    |t_1,\cdots,t_r \rangle _{\tGamma}^{(\tY, D)}=0 \text{ when } \deg \longvec{\alpha _A} +\sum_{i=1}^r{\deg t_i} > 2\vdim \moduli(Y, D, \pi_{A*} \tGamma) - 2\delta.$$

\end{corollary}
\begin{proof}
For the type I case, let $q\in H^0(Z,N_{Z/X})$ be the section found in Lemma~\ref{codimtypeI}. For the type II case, let $q$ be the zero section. Apply Lemma~\ref{topovanish} to the map $$\vp : \bM=\moduli(Bl_{q(Z)}Y, D, \tGamma) \to \moduli(Y, D, \pi_* \tGamma) \to \moduli(Y, D, \pi_{A*} \tGamma)=\bN, $$ where $U\subset \bN$ collects all relative stable maps supported away from $q(Z)$.
\begin{align*}
& \tGamma \text{ has genus-zero weight } g: V(\tGamma) \stackrel{\equiv 0}{\rightarrow} \mathbb{Z}_{\geq 0} , \text{ and } N_{Z/X} \text{ is convex}.\\
\Longrightarrow & \bM=\moduli(Bl_{q(Z)}Y, D, \tGamma) \to \moduli(Y, D, \pi_* \tGamma) \text{ has compatible perfect obstruction theories.}
\end{align*}
$\moduli(Y, D, \pi_* \tGamma) \to \moduli(Y, D, \pi_{A*} \tGamma)=\bN$ also has compatible perfect obstruction theories because it forgets $[n]-A$ legs. On the other hand, $b(v) \neq \pi^!\pi_*b(v)$ for at least one $v \in V(\tGamma)$ implies that $\vp (\bM) \cap U = \emptyset$.
The second assumption in Lemma~\ref{topovanish} follows from Lemma~\ref{codimtypeI} and Lemma~\ref{codimtypeII}.
\end{proof}

\begin{corollary}{\label{relvanish2}}
Suppose $Z \subset X$ is of type I or II. Suppose $\tGamma$ is an admissible weighted graph for $(\tY, D)$ with genus-zero weight $g: V(\tGamma) \stackrel{\equiv 0}{\rightarrow} \mathbb{Z}_{\geq 0}$. Let $\pi_{A*}$ be a composition of push-forward and the map forgetting the $[n]-A$ legs. Assume $\moduli(Y, D, \pi_{A*} \tGamma)$ makes sense.\\

If one further assumes $j\in [n]$ and $\omega_{\{j\}} \in H^*(\tY)$ with $PD_{\tY}(\omega_{\{j\}})$ sitting inside the image of  $H_*(E) \to H_*(\tY)$, where $E$ is the exceptional divisor, then we have 
$$\langle \longvec{\pi^* \alpha _A} \cdot \longvec{\tau_\bullet \gamma_{[n]}} \cdot \longvec{\omega_{\{j\}}}   |t_1,\cdots,t_r \rangle _{\tGamma}^{(\tY, D)}=0 \text{ when } \deg \longvec{\alpha _A} +\sum_{i=1}^r{\deg t_i} > 2\vdim \moduli(Y, D, \pi_{A*} \tGamma) - 2\delta.$$
\end{corollary}
\begin{proof}
Apply Corollary~\ref{vanish2} to 
$
\xymatrix{
\bM' \ar[d] \ar[r] & \bM= \moduli(\tY, D, \tGamma) \ar[d]^{\ev_j} \ar[r] & \moduli(Y, D, \pi_{A*}\tGamma)=\bN\\
E \ar[r] & \tY.
}
$\\
Any curve $[C\to \tY_l \to \tY] \in \bM'$ touches the exceptional divisor in $\tY$, therefore the correponding image in $\bN$ touches $Z \subset Y$.
Now it follows from Lemma~\ref{codimtypeI} and Lemma~\ref{codimtypeII}.
\end{proof}
\subsection{Absolute case}
Recall $\pi : \tX \to X$ is the blow up along $Z$. In the following theorem, sets $I, J, K, A$ can be empty sets. When $A$ is empty, $\deg \longvec{\alpha_A}$ will be counted as zero.
\begin{theorem}{\label{vanishingthm}}
$I, J, K$ are disjoint sets with $J \subset [n]$.
Suppose $Z = (\coprod_{i\in I} Z_i) \cup (\coprod_{j\in J} Z_j) \cup (\coprod_{k\in K} Z_k)$ is a disjoint union of submanifolds in $X$, with the following assumptions:
\begin{enumerate}
\item For each $i \in I\cup J$, $Z_i \subset X$ is either of type I or of type II.
\item For each $k \in K$, $N_{Z_k/X}$ is convex.
\item The curve class $\tbeta = \pi^!\beta + \sum_{i\in I} d_i e_i + \sum_{j\in J} d_j e_j +\sum_{k\in K} d_k e_k$ with $d_i\neq 0$ for all $i\in I$, and $0 \neq \beta\in H_2(X)$. Here $e_\bullet$ are the line classes in the corresponding exceptional divisors.
\item $\longvec{\omega_J}$ is a collection of cohomology classes in $H^*(\tX)$. And $PD_{\tX}(\omega_j)$ lies in the image of $H_*(E_j) \to H_*(\tX)$. 
\end{enumerate}
For $i\in I\cup J$, define 
$$
\delta_i= \left\{
\begin{array}{ll}
\rk(\calF) -1 & \text{ if } Z_i \subset \text{ X is of type I, and } \calF\subset N_{Z_i/X} \text{ is generated by global sections.}\\
\rk(N_{Z_i/X})-1 & \text{ if } Z_i \subset \text{ X is of type II.}
\end{array}
\right.
$$
Then 
$$\langle \longvec{\pi^*\alpha_A} \cdot \longvec{\tau_\bullet \gamma_{[n]}} \cdot \longvec{\omega_J} \rangle _{0, n, \tbeta}^{\tX}=0 \text{ when } \deg \longvec{\alpha_A} > 2\vdim \moduli_{0, A} (X, \beta)- 2\sum_{i\in I} \delta_i - 2\sum_{j\in J} \delta_j.$$
Here $\longvec{\alpha_A}$ is a collection of cohomology classes from $X$ with $A \subset [n]$, and $\longvec{\tau_\bullet \gamma_{[n]}}$ are arbitrary descendant insertions of $\tX$.
\end{theorem}

\begin{proof}
For $i \in I \cup J \cup K$, define 
$$Y_i := \bP_{Z_i} (N_{Z_i/X}\oplus \sO_{Z_i}) , \quad \pi_i: \tY_i:= Bl_{Z_i} Y_i \to Y_i , \quad D_i:= \bP_{Z_i} (N_{Z_i/X})\subset Y_i.$$
Apply the degeneration for blow-up:
$$
\xymatrix {
\tX \ar@{~}[r] \ar[d] & \tX \cup \coprod_{i\in I}\tY_i \cup  \coprod_{j\in J}\tY_j \cup \coprod_{k\in K}\tY_k   \ar[d]\\
X \ar@{~}[r] & \tX \cup \coprod_{i\in I}Y_i \cup  \coprod_{j\in J}Y_j \cup \coprod_{k\in K} Y_k
}
$$
Given $\big(\tGamma, \{ \tGamma_i\}_{i\in I} ,  \{ \tGamma_j\}_{j\in J} , \{ \tGamma_k\}_{k\in K}   \big) \in \Omega_{0, n, \tbeta}$ , we have 
$$\big(\pi_{A*}\tGamma, \{ \pi_{i,A,*}\tGamma_i\}_{i\in I} ,  \{ \pi_{j,A,*}\tGamma_j\}_{j\in J} , \{ \pi_{k,A,*}\tGamma_k\}_{k\in K}  \big) \in \Omega_{0, n, \beta},$$ 
where 
$\pi_{i,A,*}$, $\pi_{j,A,*}$ and $\pi_{k,A,*}$ are the compositions of push-forward and the map forgetting marked legs corresponding to $[n]-A$ in the absolute case. Note for each $i\in I$, $\tGamma_i$ can't be empty, and $b(\tGamma_i) = \pi_i^!b(\pi_{i,A,*}\tGamma_i) + d_i e_i$.\\

Let $\theta_0^* , \theta_i^* , \theta_j^*, \theta_k^*$ refer to the distribution of insertions to various pieces $\tX, \tY_i, \tY_j, \tY_k$. We can choose the distribution so that
\begin{enumerate} 
\item $\theta_i^*\longvec{\pi^*\alpha_A}$ , $\theta_j^*\longvec{\pi^*\alpha_A}$ and $\theta_k^*\longvec{\pi^*\alpha_A}$ are the pull back of cohomology classes from $Y_i$ , $Y_j$ and $Y_k$.
\item $\theta^* \longvec{\omega_J}$ are distributed to the corresponding divisors in $\tY_j$, for $j\in J$. 
\end{enumerate}
Argue by contradiction, suppose the invariant is not zero, then there exists
$$\big(\tGamma, \{ \tGamma_i\}_{i\in I} ,  \{ \tGamma_j\}_{j\in J} , \{ \tGamma_k\}_{k\in K}  \big) \in \Omega_{0, n, \tbeta},$$
$$\longvec{t_i}\in H^*(D_i^{\#\text{roots of } \tGamma_i}), \text{ for } i\in I\cup J \cup K,$$
so that

$$
\left\{
\begin{array}{ll}
\langle \theta_0^* \longvec{\pi^*\alpha_A} \cdot \theta_0^*\longvec{\tau_\bullet \gamma_{[n]}} | \{ \longvec{t_i} \}_{i\in I} | \{ \longvec{t_j} \}_{j\in J} | \{ \longvec{t_k} \}_{k\in K}
\rangle_{\tGamma}^{(\tX, \coprod_{i\in I\cup J \cup K} D_i)} \neq 0, \\\\
\langle \theta_i^* \longvec{\pi^*\alpha_A} \cdot \theta_i^*\longvec{\tau_\bullet \gamma_{[n]}} |  \longvec{t_i}^\vee \rangle_{\tGamma_i}^{(\tY_i, D_i)} \neq 0 \text{ for all } i\in I,\\\\
\langle \theta_j^* \longvec{\pi^*\alpha_A} \cdot \theta_j^*\longvec{\tau_\bullet \gamma_{[n]}} \cdot \theta_j^* \longvec{\omega_J} |  \longvec{t_j}^\vee \rangle_{\tGamma_j}^{(\tY_j, D_j)} \neq 0 \text{ for all } j\in J,\\\\
\langle \theta_k^* \longvec{\pi^*\alpha_A} \cdot \theta_k^*\longvec{\tau_\bullet \gamma_{[n]}} |  \longvec{t_k}^\vee \rangle_{\tGamma_k}^{(\tY_k, D_k)} \neq 0 \text{ for all } k\in K.
\end{array}
\right.
$$
Given $i\in I \cup J \cup K$, define $\deg \theta_i^* \longvec{\pi^*\alpha_A} :=\sum_{a\in A \cap \text{legs of} \tGamma_i} \deg \alpha_a$, then we have
\begin{align*}
& \deg \longvec{t_i}^\vee + \deg \theta_i^* \longvec{\pi^*\alpha_A} \leq 2\vdim \moduli (Y_i, D_i, \pi_{i,A,*}\tGamma_i)-2\delta_i
& \text{ by Corollary~\ref{relvanish1}}.\\
& \deg \longvec{t_j}^\vee + \deg \theta_j^* \longvec{\pi^*\alpha_A} \leq 2\vdim \moduli (Y_j, D_j, \pi_{j,A,*}\tGamma_j)-2\delta_j
& \text{ by Corollary~\ref{relvanish2}}.\\
& \deg \longvec{t_k}^\vee + \deg \theta_k^* \longvec{\pi^*\alpha_A} \leq 2\vdim \moduli (Y_k, D_k, \pi_{k,A,*}\tGamma_k)
& \text{ by Lemma~\ref{topovanish}}.
\end{align*}
On the other hand, by the assumption on $\deg \longvec{\alpha_A}$, we have
\begin{align*}
& \deg \theta_0^* \longvec{\pi^*\alpha_A} + \sum_{i\in I\cup J \cup K} \deg \theta_i^* \longvec{\pi^*\alpha_A} + \sum_{i\in I\cup J\cup K} \deg \longvec{t_i}^\vee + \sum_{i\in I\cup J \cup K} \deg \longvec{t_i}\\
= & \deg \longvec{\pi^*\alpha_A} + 2\sum_{i\in I\cup J \cup K} (\dim D_i)\bullet (\# \text{roots of }\tGamma_i)\\
> & 2\vdim \moduli_{0, A} (X, \beta) -2\sum_{i\in I} \delta_i -2\sum_{j\in J} \delta_j + 2\sum_{i\in I\cup J\cup K} (\dim D_i)\bullet (\# \text{roots of }\tGamma_i)\\
= & 2\sum_{i\in I\cup J\cup K} \vdim \moduli (Y_i, D_i, \pi_{i,A,*}\tGamma_i) -2 \sum_{i\in I} \delta_i -2 \sum_{j\in J} \delta_j + 2\vdim \moduli(\tX, \coprod_{i\in I\cup J\cup K} D_i, \pi_{A*}\tGamma).
\end{align*}
Combine all inequalities, we obtain
$$ \deg \theta_0^* \longvec{\pi^*\alpha_A} + \sum_{i\in I\cup J \cup K} \deg \longvec{t_i} > 2\vdim \moduli(\tX, \coprod_{i\in I\cup J \cup K} D_i, \pi_{A*}\tGamma).$$
However, $\moduli(\tX, \coprod_{i\in I\cup J\cup K} D_i, \tGamma) \to \moduli(\tX, \coprod_{i\in I\cup J\cup K} D_i, \pi_{A*}\tGamma)$ forgets $\{ \text{the marked legs of } \tGamma \} -A$, and therefore has compatible perfect obstruction theories. 
$$ \Longrightarrow \langle \theta_0^* \longvec{\pi^*\alpha_A} \cdot \theta_0^*\longvec{\tau_\bullet \gamma_{[n]}} | \{ \longvec{t_i} \}_{i\in I} | \{ \longvec{t_j} \}_{j\in J} | \{ \longvec{t_k} \}_{k\in K}
\rangle_{\tGamma}^{(\tX, \coprod_{i\in I\cup J \cup K} D_i)} = 0 \text{ by Lemma~\ref{topovanish}},$$
which is a contradiction.
\end{proof}
\begin{example}{\label{surface}}
Suppose $X$ is an algebraic surface, which is neither rational nor ruled. Let $X_0$ be the minimal model of $X$. Since GW-invariants are deformation invariant, we may assume $\pi: X \to X_0$ is the blow-up at $r$ distinct points $a_1, \cdots, a_r$. Suppose $0 \neq \beta \in H_2(X_0)$.
$$K_{X_0} \text{ is nef } \Longrightarrow \vdim\moduli_{0,0}(X_0, \beta)= (2-3) +0 - \beta\cap K_{X_0} \leq -1.$$
Assume $\tbeta = \pi^!\beta+ \sum_{k=1}^r d_k e_k$, where $d_k \in \mathbb{Z}$. We apply the previous theorem to the case $Z=\coprod_{k\in K} Z_k= \{ a_1, \cdots, a_r\}$, with the set $I=J=A=\emptyset$. We have 
$$\deg \longvec{\alpha_A}=0 > -2 \geq  2\vdim\moduli_{0,0}(X_0, \beta).$$
By the previous theorem, $g=0$ descendant GW-invariants of $X$ are all zero if $\beta\neq 0$.
Since exceptional divisors are disjoint,
$$
\langle \tau_{a_1}\gamma_1, \cdots, \tau_{a_n}\gamma_n \rangle _{0,n, \tbeta}^X=
\left\{
\begin{array}{ll}
\text{ invariants around the exceptional divisor }\bP^1 & ,\text{ if } \tbeta=d_k e_k \text{ for some } k, \text{ with } d_k>0.\\
0 & ,\text{ otherwise.}
\end{array}
\right.
$$
The first case can be computed by obstruction bundles.\\
When $p_g(X)>0$, this result can also be deduced from Image Localization Theorem in ~\cite{structure} (see also ~\cite{thetaloc} for algebro-geometric analogue)in symplectic geometry. In fact, Image Localization Theorem is much more powerful than our argument because it can also handle higher genus GW-invariants when $p_g>0$.
\end{example}

\bigskip

\begin{example}
Suppose $K_X$ is nef, and $Z$ is a smooth curve in $X$ with genus $g(Z) \geq 1$. Then we have zero descendant GW-invariants 
$$\langle \cdots \rangle_{0,n, \tbeta}^{\tX} \equiv 0 \text{, when } \tbeta= \pi^!\beta + d e \in H_2(\tX) \text{ with } \beta\neq 0 \text{ and } d \neq 0.$$
To see this, note $\vdim \moduli_{0,0}(X, \beta) = \dim X -3+0 - \beta \cap K_X \leq \dim X -3$. Apply the vanishing theorem to $J=K=A =\emptyset$, then $\delta = \codim (Z,X)-1 = \dim X -2$. Hence $\deg \longvec{\alpha_A}=0 > -2 \geq 2\vdim \moduli_{0,0}(X, \beta) - 2\delta$.
\end{example}

\bigskip

\begin{example}{\label{counter2}}
Let $Z=\bP^2$, and $X$ is the projective completion of $\sO(-3)\oplus \sO(-3) \to Z$. This example shows Theorem~\ref{vanishingthm} doen't not hold for arbitrary blow-ups. Let $\pi : \tX \to X$ be the blow-up along $Z$. The exceptional divisor is $E\cong Z \times \bP^1$ with normal bundle $N_{E/\tX} \cong \sO_Z(-3) \boxtimes \sO(-1)$. Let $[\ell_1]$ and $[\ell_2]$ be the line classes in $Z$ and $\bP^1$. 
$$E \stackrel{i}{\rightarrow} \tX \stackrel{\pi}{\rightarrow} X \stackrel{p}{\rightarrow} Z.$$
Then $\pi^![\ell_1]= i_*(\ell_1 -3 \ell_2)$. Consider $\moduli_{0,1}(\tX, i_*(d\ell_1)) \to \moduli_{0,0}(X, d\ell_1)$ with $d\geq 1$. Let $I=J=A=\emptyset$ in Theorem~\ref{vanishingthm}. We have $\deg \longvec{\alpha_A} = 0 > 2(1-3d) =2\vdim \moduli_{0,0}(X, d\ell_1)$. If Theorem~\ref{vanishingthm} holds in this example, then it implies all GW-invariants of $\moduli_{0,1}(\tX, i_*(d\ell_1))$ are zero.\\
On the other hand,
$(E + 3\pi^*p^*H_1)|_E= -H_2 \in H^2(E)$, where $H_1$ and $H_2$ are hyperplane classes of $Z$ and $\bP^1$ in $E$. Let $\calU_d \to \moduli_{0,0}(Z, d\ell_1)$ be the obstruction bundle associated to $\sO(-3) \to Z$.
\begin{align*}
\langle 3\pi^*p^*H_1 \wedge (E + 3\pi^*p^*H_1) \rangle_{0,1, d\ell_1}^{\tX}= & -\int_{\moduli_{0,1}(E, d\ell_1)^\vir} c_{top}(\calU_d)\cap \ev^* (H_1 \boxtimes H_2)\\
= & -d \int_{\moduli_{0,0}(Z, d\ell_1)^\vir} c_{top}(\calU_d) = -d \cdot K_d.
\end{align*}
The number $K_d$ has been computed in~\cite{mirror1}, and is non-zero in general (e.g. $K_1=3$).
\end{example}

\bigskip

\begin{example}
Suppose $N_{Z/X}$ is generated by global sections and has rank $r$. Let $E$ be the exceptional divisor of $\pi: \tX \to X$.
Given $a_i\geq0$, $0 \neq \beta\in H_2(X)$ and $\alpha_i \in H^*(X)$, then 
$$\langle E^{a_1}\pi^*\alpha_1, E^{a_2}\pi^*\alpha_2, \cdots, E^{a_n}\pi^*\alpha_n \rangle_{0,n, \pi^!\beta}^{\tX}=0 \text{ when } 0<\sum_{i=1}^n a_i < r-1.$$
To see this, may assume $a_1 > 0$, then apply Theorem~\ref{vanishingthm} to:
$$
\left\{
\begin{array}{l}
K =\emptyset \text{ and } \delta=r-1,\\
J=\{ 1 \}\subset [n] \text{ with } \omega_1=E^{a_1},\\
\longvec{\alpha_A}=(\alpha_1, \cdots, \alpha_n) \text{ with } A=[n]\\
\longvec{\gamma_{[n]}}= \text{ all remaining insertions}.
\end{array}
\right.
$$
Then $\displaystyle \deg \longvec{\alpha_A} = \vdim\moduli_{0,n}(\tX, \pi^!\beta) - \sum_{i=1}^n a_i > \vdim \moduli_{0,n}(X, \beta) -(r-1)$. Therefore the invariant vanishes. One can use the similar argument to $\moduli_{0,n}(\tX, \pi^!\beta) \to \moduli_{0,n-m}(X, \beta)$ and show that if $1\leq m \leq n $, then
$$\langle E^{a_1}, E^{a_2}, \cdots, E^{a_m}, E^{a_{m+1}}\pi^*\alpha_{m+1}, \cdots \cdots, E^{a_n}\pi^*\alpha_n \rangle_{0,n, \pi^!\beta}^{\tX}=0 \text{ when } \sum_{i=1}^n a_i < r-1 + m .$$

If there are too many insertions coming from the exceptional divisor, then the invariant may not vanish.
For example, take $X=\bP^3$ and $Z=$ a point. Let $[\ell]$ be a line class in $X$. A computation in~\cite{Gathmann0} shows 
$$\langle E^2, E^2, \cdots, E^2 \rangle_{0,12, 3\ell}^{\tX} = -2332 \neq 0.$$
\end{example}

\begin{remark}
\mbox{}
\begin{enumerate}
\item
Suppose $Z=Z_1\coprod Z_2$ with $N_{Z_i/X}$ both generated by global sections. Let $\pi : \tX= Bl_Z X \to X$ and $\pi_i : \tX \to X_i=Bl_{Z_i}X$. To test if a GW-invariant of $\tX$ vanishes or not, using different base manifolds can yield different vanishing criteria. For example, let $\tbeta=\pi^!\beta + d_1 e_1 + d_2 e_2 \in H_2(\tX)$ with $d_1, d_2 > 0$.
If $\alpha_a\in H^*(X_1)$, then
$$\langle \longvec{\pi_1^*\alpha_A} \cdot \longvec{\tau_\bullet \gamma_{[n]}} \rangle _{0, n, \tbeta}^{\tX}=0 \text{ when } \deg \longvec{\alpha_A} > 2\vdim \moduli_{0, A} (X_1, \beta + d_1 e_1)- 2\rk(N_{Z_2/X}).$$
However, this result can not be deduced from the vanishing criterion for $\tX \to X$ because $\longvec{\alpha_A}$ may not come from cohomology classes of $X$.
\item
It is not necessary to test all possible base manifolds. In Theorem~\ref{vanishingthm}, suppose $I= I_+ \coprod I_-$ such that 
$\displaystyle \tbeta = \pi^!\beta + \sum_{i\in I} d_i e_i + \sum_{j\in J} d_j e_j +\sum_{k\in K} d_k e_k$ with $d_i > 0$ for all $i\in I_+$, and $d_i <0$ for all $i\in I_-$.

A simple argument shows : if an invariant of $\tX$ satisfies the vanishing criterion for $\tX \to X$, then it automatically satisfies the vanishing criterion for $\displaystyle \tX \to Bl_{(Z_{I_-})} X$, where $Z_{I_-} =\displaystyle \coprod_{i \in I_-} Z_i$.
\end{enumerate}
\end{remark}

\bibliographystyle{plain}
\bibliography{GW}
Department of Mathematics, 
Brandeis University,
415 South Street
MS 050,
Waltham, MA 02454\\
\textit{E-mail address}: hhlai@brandeis.edu
\end{document}